\documentclass{article}

\usepackage{amssymb,amsmath,eufrak,amsthm,amscd}
\usepackage[matrix,arrow,curve]{xy}

\usepackage{hyperref}

\newtheorem{theorem}{Theorem}
\newtheorem{lemma}[theorem]{Lemma}
\newtheorem{proposition}[theorem]{Proposition}
\newtheorem{corollary}[theorem]{Corollary}

\newtheorem*{statement*}{Statement}
\newtheorem*{theorem*}{Theorem}
\newtheorem*{lemma*}{Lemma}
\newtheorem*{fact*}{Fact}

\newcommand{\diffspec}{\operatorname{Spec^\Delta}}

\newcommand{\diffsmax}{\operatorname{SMax^\Delta}}

\newcommand{\Qt}{\operatorname{Qt}}

\newcommand{\diffaut}{\operatorname{Aut^\Delta}}
\newcommand{\diffgal}{\operatorname{Gal^\Delta}}

\theoremstyle{definition}
\newtheorem{definition}[theorem]{Definition}
\newtheorem*{definition*}{Definition}
\newtheorem{example}[theorem]{Example}

\newtheorem*{example*}{Example}

\theoremstyle{remark}

\newtheorem*{remark*}{Remark}

\author{D.\,V.~Trushin}
\title{Splitting fields and general differential Galois theory\footnote{\Xy-pic package is used}}
\date{}

\begin{document}

\maketitle

\begin{abstract}

An algebraic technique is presented that does not use results of
model theory and makes it possible to construct a general Galois
theory of arbitrary nonlinear systems of partial differential
equations. The algebraic technique is based on the search for prime
differential ideals of special form in tensor products of
differential rings. The main results demonstrating the work of the
technique obtained are the theorem on the constructedness of the
differential closure and the general theorem on the Galois
correspondence for normal extensions.

\end{abstract}

\section{Introduction}

The first algebraic proof of the existence and uniqueness of a
differential closure for a differential field with finitely many
pairwise commuting derivations appeared in Kolchin's
paper\footnote{Kolchin uses the term constrained
closure}~\cite{KL3}. A fairly detailed overview of papers on
differentially closed fields is contained in the introduction of
that paper; however certain questions on the structure of the
differential closure remained unanswered. On the other hand, back in
Shelah's paper~\cite{Sh} the method related to complete totally
transcendental first-order theories was used, which contains
somewhat more information about the differential closure. In
Kolchin's opinion (see~\cite{KL3}, p.~141), Shelah's proof is one of
the most difficult. Kolchin's paper~\cite{KL3} is an attempt to
carry out all the proofs in the language of differential algebra.
But we are going to expound Shelah's proof in the case of
differential fields of characteristic zero without using the
model-theoretic technique. For that we develop a fairly simple
technique of searching for differential ideals in tensor products of
differential rings. Thus, we not only translate into the algebraic
language the fairly powerful proof of Shelah, but also make a
substantial simplification of it.

The model-theoretic  proof is divided into two stages; we shall
discuss it only in respect to the theory of differential fields with
finitely many derivations. At first the motion of a constructed
models is introduced and Ressayre's theorem on the uniqueness of a
constructed model is proved for complete theories\footnote{we shall
use Poizat's book instead of the original paper of Shelah}
(see~\cite{Pz}, Theorem~10.18). Then for totally transcendental
theories all simple models are shown to be constructed
(see~\cite{Pz}, \textsection\textsection~18.1, 18.2). In the present
paper we follow the same scheme of proof. The strongest result in
this direction is the proof of Theorem~\ref{DCFuniquetheorem} on the
constructedness of the differential closure. Moreover, the algebraic
technique obtained makes it possible to advance further and to
obtain the notion of a splitting field for an arbitrary system of
differential equations, and to prove for splitting fields the
theorem of the Galois correspondence for normal extensions.

At present there exist a vast number of papers on differential
Galois theory of various types of differential equations. For
example, the Picard-Vessiot theory of linear differential equations
(see~\cite{SW}) or the Picard-Vessiot theory of parametrized linear
differential equations (see~\cite{SC}). One of the drawbacks of
these papers is the superfluous assumption about the properties of
the initial differential field; namely, in the Picard-Vessiot theory
it is assumed that the constants are algebraically closed, and in
the theory of parametrized equations even more -- the differential
closedness. In applications of these theories such assumptions are
inadmissible, especially the assumption of differential closedness,
since the initial fields usually are fields constructed of function,
and their differential closure are too big and have fairly
complicated structure. the notion of a splitting field that we
introduce removes this defect; namely, an ordinary Picard-Vessiot
extension (and its parametrized analogue) is a special case of a
splitting field, for the construction of which no additional
assumptions on the initial field are required.

A detailed survey the existing Galois theories can be found in the
paper of Cassidy and Singer (see~\cite{SC}, p.~145-137). However, it
seems that the most interesting is the cycle of papers of
Pillay~\cite{Pl1,Pl2,Pl3}, in which he completely describes the
structure of the differential closure in terms of the algebraic
closure, generalized Galois extensions, and extensions by trivial
elements corresponding to symmetric groups. In the end it is
desirable to find a purely algebraic proof of Pillay's results. The
present paper is merely the first step on this road. However,
already here we succeeded in constructing an algebraic apparatus
that makes it possible to avoid using the model-theoretic technique
ar certain stages. Furthermore, it is worth mentioning that the
technique developed by us finds an unexpected application in
difference algebra; namely, one can use this technique for proving
the existence of difference-closed fields in a fairly board
sense\footnote{http://arxiv.org/abs/0908.3865}.

the concluding detail in the present paper is Theorem~\ref{Variety},
in which we discuss the connection of the differential spectrum of a
differentially finitely generated algebra over a field with the
corresponding differentially algebraic variety. It is shown that
from this viewpoint the set of locally closed points of the
differential spectrum plays the same role as the maximal spectrum
for finitely generated algebras over a field in commutative case.

This research was initiated by A.~I.~Ovchinnikov. He acquainted the
author with the book of Poizat~\cite{Pz} and posed the questions
that are solved in this paper. In the preparation of the paper
invaluable help was given by M.S.~Singer. Due to his comments, a
detailed example emerged that illustrates effects of various kinds.
I am grateful to my supervisor A.V.~Mikhalev.

\section{Detailed structure of the paper}

The paper consists of three key parts. The first is
\textsection\textsection~\ref{sec3} and~\ref{sec4}; here all the
main technique is collected on which the further exposition is
based. We consider questions of constructing differential ideals
with requisite properties, and on the basis of them --- constructing
differential closure of special types. Sections~\ref{sec5}
and~\ref{sec6} are devoted to the translation of model-theoretic
results into the algebraic language; these sections constitute the
second part of the text. Subsequently we define the notion of a
splitting field of an arbitrary system of differential equations and
the notion of a differential Galois group.
Sections~\ref{sec7}-\ref{sec9} are devoted to the study of the new
notions and their connection with the known theories.

We now describe the structure of the paper in more details.
Section~\ref{sec3} is devoted to the terminology used.
Section~\ref{sec4} is the a `core' of the technique that is
developed. It consists of three subsections, in each of which we
construct specific prime differential ideals in tensor products of
differential rings. Subsection~\ref{sec41} is devoted to
construction a prime differential ideals with residue field that
does not contain new constants (theorem~\ref{algconst}), and this
technique is applied to the Picard-Vessiot theory
(Proposition~\ref{PicVes}) and to constructing a differential
closure with algebraic field of constants
(Theorem~\ref{DCFalgconst}). In Subsection~\ref{sec42} we consider
prime differential ideals with a slightly more subtle property
(Theorem~\ref{atomic}). As in Subsection~\ref{sec41}, based on the
properties of the ideals obtained, we construct the corresponding
differential closure (Theorem~\ref{DCFatomic}).
Subsection~\ref{sec43} is devoted to prime differential ideals that
have a certain universal property (Theorem~\ref{univers}). Such
ideals enable us to construct constructed fields and constructed
differential closures (Theorem~\ref{DCFunivers}).

Section~\ref{sec5} is devoted to the detailed study of constructed
fields. Our main goal is the translation into the algebraic language
of the model-theoretic proof of Ressayre's theorem --- Theorem~10.18
in~Ch.~10 of~\cite{Pz} (see~Theorem~\ref{ConstructiveUnique}). In
\textsection~\ref{sec6} we translate into the algebraic language the
proof of Theorem~18.1 in~Ch.~18 of~\cite{Pz}
(see~Theorem~\ref{DCFuniquetheorem}).

Due to the results of \textsection\textsection~\ref{sec5}
and~\ref{sec6} we are able to define in~\textsection~\ref{sec7} the
notion of a splitting field for an arbitrary system of differential
equations. In Subsection~\ref{sec71} we define an abstract splitting
field, study in detail its properties, and show the connections with
already known notions. Propositions~\ref{ExistSplitClose}
and~\ref{ExtAutoSplit} are devoted to the existence and uniqueness.
Subsection~\ref{sec72} us devoted to the connections of splitting
fields with one another. All this leads to the notion of a normal
extension, which is presented in Subsection~\ref{sec73} and
constitutes the main object of study in Galois theory.
In~\textsection~\ref{sec8} the theorem on the Galois correspondence
is proved for normal extensions of differential fields
(Theorem~\ref{GaloisCorr}) and an example is given, which
graphically illustrates the behavior of the objects that we define.

In~\textsection~\ref{sec9} we give a geometric interpretation of
locally closed points of the differential spectrum of a
differentially finitely generated algebra over an arbitrary
differential field (Theorem~\ref{Variety}) in terms of the
corresponding differential algebraic varieties and the Galois group
of the differential closure.

\section{Definitions and notation}\label{sec3}

Throughout the paper we consider some differential ground field $K$
of characteristic zero. All differential rings are assumed to be
algebras over $K$ with a finite set of pairwise commuting
derivations
$$
\Delta=\{\,\delta_1,\ldots,\delta_m\,\}.
$$
{\it Differentially finitely generated algebras} are by definition
differentially finitely generated algebras over $K$. A {\it simple
ring} is by definition a differential ring without nontrivial
differential ideals. If it is not specified over which ring a tensor
product is taken, then we mean the tensor product over the field
$K$. The differential spectrum of a ring $A$ is denoted by
$\diffspec A$, and $\diffsmax A$ denotes the set of locally closed
points of the differential spectrum. The field $(A/\frak p)_\frak p$
is called the {\it residue field of a prime ideal} $\frak p\subseteq
A$. The field of fractions of an integral domain $B$ is denoted by
$\Qt(B)$. If in a ring $D$ two subrings $A$ and $B$ are given, then
$A\cdot B$ denotes the smallest subring of $D$ generated by $A$ and
$B$. The definition and notation in commutative algebra are the same
as in~\cite{AM}, and in differential algebra as in~\cite{Kl}.

\section{Ideals of tensor products}\label{sec4}

In this section we construct special types of differential ideals in
tensor products of differential algebras. Search for prime
differential ideals of special types means search for composites of
differential fields of special type. Such composites appear, in
particular, in the construction of differentially closed fields
containing some given differential field. Therefore, the more
sophisticated ideals we can find, the more delicate differential
closed overfield we can construct. Here we search for differential
prime ideals of three types. By using each of them we show the
existence of special types of differentially closed fields. We begin
with the most rough condition -- the absence of new constants, which
works well for Picard-Vessiot extensions, and conclude with
constructing a differential closure of a fixed field.

Let $\{\,B_\alpha\,\}_{\alpha\in\Lambda}$ be some set of simple
$K$-algebras. Suppose that in addition the algebras $B_\alpha$ are
differentially finitely generated. We consider  $R=\otimes_\alpha
B_\alpha$ as the inductive (direct) limit of finite tensor products
over $K$ (see~\cite[Ch.~2, Exercise~23]{AM}).

\subsection{Algebraic constants}\label{sec41}

To begin with, we isolate the simplest class of differential ideals
of $R$ related to constants.

\begin{theorem}\label{algconst}
Let $\{B_\alpha\}_{\alpha\in \Lambda}$ be some set of simple
differentially finitely generated $K$-algebras. We denote by $C$ the
subfield of constants of the field $K$. Then there exists a prime
differential ideals $\frak p$ of $R=\otimes_\alpha B_\alpha$ such
that the constants of the residue field of $\frak p$ are algebraic
over $C$.
\end{theorem}
\begin{proof}

Consider the set $\mathbb S$ consisting of pairs
$(\otimes_{\theta\in\Theta}B_\theta,\frak p_\Theta)$, where $\Theta
\subseteq \Lambda$ and $\frak p_\Theta$ is a prime differential
ideal in $\otimes_{\theta\in\Theta}B_\theta$ such that the constants
of its residue field are algebraic over $C$. We order the set
$\mathbb S$ as follows:
$$
(\otimes_{\theta\in\Theta_1}B_\theta,\frak p_{\Theta_1})\leqslant
(\otimes_{\theta\in\Theta_2}B_\theta,\frak p_{\Theta_2})
\Leftrightarrow \Theta_1\subseteq\Theta_2, \, \frak
p_{\Theta_2}\cap\otimes_{\theta\in\Theta_1}B_\theta=\frak
p_{\Theta_1}
$$
Note that $\mathbb S$ is not empty; for example, it contains the
pair $(B_\alpha, 0)$ for any elements $\alpha$ (see~\cite[Ch.~III,
\textsection~10, Proposition~7(d)]{Kl}). The set thus constructed
satisfies the hypothesis of Zorn's lemma, for which one must take
the direct limit over a chain. Let
$(\otimes_{\theta\in\widehat\Theta}B_\theta,\frak
p_{\widehat\Theta})$ be a maximal element of the set $\mathbb S$. We
claim that $\otimes_{\theta\in\widehat\Theta}B_\theta$ coincides
with the whole ring $R$.

We introduce the notation
$$
R_{\widehat\Theta}=\otimes_{\theta\in\widehat\Theta}B_\theta/\frak
p_{\widehat\Theta},\,\,
B_{\widehat\Theta}=\otimes_{\theta\in\widehat\Theta}B_\theta,\,\,S=B_{\widehat\Theta}\setminus\frak
p_{\widehat \Theta}.
$$
Consider the ring $R'=B_{\widehat \Theta}\otimes B_\alpha$. We
define $\frak p'$ to be a maximal differential ideal of it
contracting to $\frak p_{\widehat\Theta}$. We shall show that
$(R',\frak p')$ is contained in $\mathbb S$; then we shall obtain a
contradiction with the maximality. Indeed, $S^{-1}R'/\frak p'$ is a
simple differentially finitely generated algebra over the field
$S^{-1}R_{\widehat\Theta}$, the constants of which are algebraic
over the field of constants of the field $K$ (see~\cite[Ch.~III,
\textsection~10, Proposition~7(d)]{Kl}).  Therefore the constants of
the residue field of $\frak p'$ are algebraic over $C$.
\end{proof}

In the following definition we assume that a single derivation acts
on differential rings.

\begin{definition}\label{defUPV}
A differential algebra $D$ over the field $K$ is called a {\it
universal Picard-Vessiot ring} if the following properties hold:
\begin{enumerate}
\item $D$ is a simple differential ring;
\item for any linear differential equation of the form $y'=Ay$,
where $A$ is a matrix of elements of the field $K$, there exists a
fundamental matrix $F$ of solutions with coefficients in $D$;
\item as a $K$-algebra, $D$ is generated by the elements of all
fundamental matrices $F$ and $\det(F)^{-1}$.
\end{enumerate}

This definition can be found in \textsection~1 of Ch.~10
in~\cite{SW}.
\end{definition}

The following statement is a simple corollary of the previous fact.
(Definition of universal Picard-Vessiot extension is
in~\cite[chapter~10, sec.~1]{SW}.)

\begin{proposition}\label{PicVes}
For any ordinary differential field $K$ with algebraically closed
subfield of constants, there exists a universal Picard-Vessiot
extension.
\end{proposition}
\begin{proof}
Let $\{\,B_\alpha\,\}$ be the family of all Picard-Vessiot
extensions up to isomorphism. In the ring $R=\otimes_\alpha
B_\alpha$ consider the ideal $\frak p$ from theorem~\ref{algconst}.
Then the ring $R/\frak p$ has all the required properties of
Definition~\ref{defUPV}, except for, possible, 1). However, since
$R/\frak p$ does not contain new constants, it follows from
Definition~1.21 and Proposition~1.22 in \textsection~3 of~Ch.~1
in~\cite{SW} that $R/\frak p$ is the direct limit of simple
differential rings and, consequently, is itself also simple.
\end{proof}

\begin{definition}\label{def_D_closed}
A field is said to be differentially closed if any simple
differentially finitely generated algebra over this field coincides
with it.
\end{definition}

This definition can be found in~\cite{KL3}, \cite{SC},
Definition~3.2, \cite{JK}, p.~4490, \cite{Scan}, \textsection~4,
p.~28, or~\cite{MG}. A slightly more nontrivial application of
Theorem~\ref{algconst} is the following theorem.

\begin{theorem}\label{DCFalgconst}
Let $K$ be a differential field with subfield of constants $C$. Then
there exists a differentially closed field $L$ containing $K$ such
that the field of constants of $L$ is the algebraic closure of the
field $C$.
\end{theorem}
\begin{proof}
Consider the set $\{\,B_\alpha\,\}$ of all simple differentially
finitely generated algebras over $K$ up to isomorphism. Consider an
ideal $\frak p$ as in Theorem~\ref{algconst}. Let $L_1$ be the
residue field of the ideal $\frak p$. By performing the same
procedure for $L_1$, we obtain $L_2$, for $L_2$ we obtain $L_3$, and
so on. As a result we obtain the ascending chain of differential
fields
$$
K=L_0\subseteq L_1\subseteq\ldots \subseteq L_n\subseteq\ldots
$$
We set $L=\cup_k L_k$. We see that the field thus constructed does
not contain non-algebraic constants. Note that the field $L_{k+1}$
is constructed so that for any nontrivial differential ideal of the
ring of differential polynomials over the field $L_k$ there exists a
common zero in the field $L_{k+1}$.

We claim that $L$ is differentially closed. Let
$$
B= L\{y_1,\ldots,y_n\}/\frak m
$$
be a simple differentially finitely generated algebra. The
Ritt-Raudenbush basis theorem says that there exists a finite set of
differential polynomials $F$ such that $\frak m = \{F\}$. But the
set of all coefficients of elements of $F$ is finite and therefore
all of them are defined over some $L_k$ for a suitable index $k$. By
the construction of $L_{k+1}$ the family $F$ has a common zero
$\overline{a}\in L_{k+1}\subseteq L$ and therefore a differential
homomorphism $B\to L$ is defined by the rule $y_i\mapsto a_i$. Since
$B$ is simple and contains $L$ as a subalgebra, the latter
homomorphism is an isomorphism.
\end{proof}

\subsection{Local simplicity}\label{sec42}

First of all we introduce one definition.

\begin{definition}\label{def_locsim}
We say that a differential ring $B$ is {\it locally simple} if there
exists an element $s\in B$ such that the ring $B_s$ is a simple
differential ring. A prime differential ideal $\frak p$ of a
differential ring $B$ is called {\it locally maximal} if the algebra
$B/\frak p$ is locally simple.

\end{definition}

Such ideals are precisely the locally closed points of the
differential spectrum. In the terminology of Kolchin such algebras
are called {\it constrained} (see~\cite{KL3}).

To obtain more delicate results we need a more subtle technique. The
main purpose of this subsection is to prove theorem~\ref{atomic}.
The following lemma contains all the main technique.

\begin{lemma}[``on splitting'']\label{split}
Suppose that $D$ is a differential $K$-algebra that is an integral
domain, and $A$ and $B$ are differential $K$-subalgebras of $D$ such
that $D=A\cdot B$ and $B$ is differentially finitely generated. Then
there exist a differentially finitely generated $K$-subalgebra $C$
of $A$ and an element $s\in C\cdot B$ such that
$$
D_s=A\otimes_{C}(C\cdot B)_s.
$$
\end{lemma}
\begin{proof}

We choose differential generators $b_1,\ldots,b_n$ in the algebra
$B$. Then the algebra $D$ is generated over $A$ by the chosen
elements as a differential ring.

Consequently, $D$ has the form  $A\{y_1,\ldots,y_n\}/\frak p$, where
$\frak p$ is some prime differential ideal. Let
$G=\{\,f_1,\ldots,f_m\,\}$ be a characteristic set of the ideal
$\frak p$ for some fixed ranking and let $h$ be the product of the
initials and separants of elements of $G$. Thus we have
$$
D_h=A\{y_1,\ldots,y_n\}_h/[f_1,\ldots,f_n].
$$
Elements of $G$ depend on finitely many coefficients; we choose the
differential subalgebra $C$ generated by them. Then, by tensor
multiplying the exact triple
$$
0\rightarrow[f_1,\ldots,f_m]\rightarrow
C\{y_1,\ldots,y_n\}_h\rightarrow
C\{y_1,\ldots,y_n\}_h/[f_1,\ldots,f_n]\rightarrow0
$$
by $A$ over $C$ we obtain the required splitting.
\end{proof}

From the splitting lemma we immediately obtain the following
corollary.

\begin{corollary}\label{splitting_corollary}
Suppose that  $D$ is a differential $K$-algebra that is an integral
domain, and $A$ and $B$ are differential $K$-subalgebras of it such
that $D=A\cdot B$ and $B$ is differentially finitely generated. Then
for any element $h\in D$ there exist a differentially finitely
generated $K$-subalgebra $C$ in $A$ and an element $s\in C\cdot B$
such that $h\in C\cdot B$ and
$$
(A\cdot B)_{sh}=A\otimes_{C}(C\cdot B)_{sh}.
$$
\end{corollary}

We see that $C$ can be replaced by any larger differential ring.

\begin{corollary}\label{splitting_field_corollary}
Suppose that $D$ is a differential $K$-algebra that is an integral
domain, and $A$ and $F$ are differential $K$-subalgebras of it such
that $D=A\cdot F$ and $F$ is field that is differentially finitely
generated over $K$. Then for any element $h\in D$ there exist a
differentially finitely generated $K$-algebra $C$ in $A$ and an
element $s\in C\cdot F$ such that $h\in C\cdot F$ and
$$
(A\cdot F)_{sh}=A\otimes_{C}(C\cdot F)_{sh}.
$$
\end{corollary}
\begin{proof}

It follows from the hypothesis that there exists a differentially
finitely generated $K$-algebra such that $F$ is its field of
fractions. By applying the preceding corollary we obtain
$$
(A\cdot B)_{sh}=A\otimes_{C}(C\cdot B)_{sh}.
$$
If $S=B\setminus\{\,0\,\}$, then
\begin{multline*}
(A\cdot F)_{sh}=(A\cdot S^{-1}B)_{sh}=S^{-1}((A\cdot
B)_{sh})=\\
=S^{-1}(A\otimes_{C}(C\cdot B)_{sh})=A\otimes_{C}(C\cdot
S^{-1}B)_{sh}=A\otimes_{C}(C\cdot F)_{sh}.
\end{multline*}
\end{proof}

\begin{proposition}\label{surjstate}
Let $A$ and $B$ be a differential $C$-algebras, where $C$  is a
simple differential ring. Then the canonical map
$$
\diffspec A\otimes_C B\to\diffspec A \times\diffspec B
$$
is surjective.
\end{proposition}
\begin{proof}
Let $S$ be the set of nonzero elements of $C$. Since $C$ is simple,
then for any differential $C$-algebra $D$ we have the equality
$$
\diffspec D =\diffspec S^{-1}D.
$$
Thus, it suffices to prove the assertion for the case where $C$ is a
field. Let $(\frak p,\frak q)$ be an arbitrary pair of prime
differential ideals in the product indicated, and let  $S=A\setminus
\frak p$ and $T=B\setminus \frak q$. Then it follows from
Exercise~21 in~Ch.~3 of~\cite{AM} and the properties of tensor
product that the inverse image of this pair under that map is
naturally homeomorphic to
$$
\diffspec S^{-1}(A/\frak p)\otimes_C T^{-1}(B/\frak q)
$$
Since $C$ is a field the latter ring is nonzero and therefore its
differential spectrum is not empty.
\end{proof}

\begin{theorem}\label{atomic}
Let $\{B_\alpha\}_{\alpha\in \Lambda}$ be some set of simple
differentially finitely generated $K$-algebras. Then there exists a
prime differential ideal $\frak p$ of the ring $R=\otimes_\alpha
B_\alpha$ such that any differentially finitely generated
$K$-subalgebra $B$ in the residue field of the ideal $\frak p$ is
locally simple.
\end{theorem}
\begin{proof}

Similarly to the proof of Theorem~\ref{algconst} we consider the set
$\mathbb S$ consisting of pairs
$(\otimes_{\theta\in\Theta}B_\theta,\frak p_\Theta)$, where $\Theta
\subseteq \Lambda$, $\frak p_\Theta$ is a prime differential ideal
of $\otimes_{\theta\in\Theta}B_\theta$ such that every
differentially finitely generated subalgebra $B$ in the residue
field of the ideal $\frak p_\Theta$ is locally simple. We order this
set as follows:
$$
(\otimes_{\theta\in\Theta_1}B_\theta,\frak p_{\Theta_1})\leqslant
(\otimes_{\theta\in\Theta_2}B_\theta,\frak p_{\Theta_2})
\Leftrightarrow \Theta_1\subseteq\Theta_2, \, \frak
p_{\Theta_2}\cap\otimes_{\theta\in\Theta_1}B_\theta=\frak
p_{\Theta_1}
$$
Note that $\mathbb S$ is not empty; for example, it contains
$(B_\alpha,0)$ for any element $\alpha$ (see, for
example~\cite[Ch.~III, \textsection~10, Proposition~7(b)]{Kl}). The
set thus constructed satisfies the hypothesis of Zorn's lemma. Let
$(\otimes_{\theta\in\widehat\Theta}B_\theta,\frak
p_{\widehat\Theta})$ in $\mathbb S$ be a maximal element of it. We
claim that $\otimes_{\theta\in\widehat\Theta}B_\theta$ coincides
with the whole ring $R$. We introduce the notation
$$
R_{\widehat\Theta}=\otimes_{\theta\in\widehat\Theta}B_\theta/\frak
p_{\widehat\Theta},\,\,
B_{\widehat\Theta}=\otimes_{\theta\in\widehat\Theta}B_\theta,\,\,
S=B_{\widehat\Theta}\setminus\frak p_{\widehat \Theta},\,\,
A=S^{-1}B_{\widehat\Theta}.
$$
Consider the ring $R'=B_{\widehat \Theta}\otimes B_\alpha$. We
define  $\frak p'$ as a maximal differential ideal of this ring
contracting to  $\frak p_{\widehat\Theta}$. We shall show that
$(R',\frak p')$ is contained in $\mathbb S$, then we shall obtain a
contradiction with the maximality.

We denote the ring $S^{-1}(R'/\frak p')$ by $D$, and let $B$ be an
arbitrary differentially finitely generated $K$-subalgebra of the
residue field of the ideal $\frak p'$. Since $B$ is differentially
finitely generated, there exists $h\in D$ such that $B\subseteq
D_h$. We set $B'=B\cdot B_\alpha$. By construction we have the
equation
$$
D_h=(S^{-1}A\cdot B')_h.
$$
By applying corollary~\ref{splitting_corollary} we obtain
$$
D_{sh}=S^{-1}A\otimes_{C}(C\cdot B')_{sh},
$$
for suitable $s$ and $C$. But $C$ is a differentially finitely
generated $K$-algebra in $S^{-1}A$ and therefore $C$ is locally
simple. By inverting one element, we can assume that $C$  is simple.
Then Proposition~\ref{surjstate} guaranties that the differentially
finitely generated $K$-algebra $(C\cdot B')_{sh}$  is simple. Then
it follows from Proposition~7(b) in \textsection~10 of Ch.~III
in~\cite{Kl} that $B$ is locally simple. We have obtained a
contradiction with the maximality.
\end{proof}

It should be noted that the ideal constructed in
theorem~\ref{atomic} is a special case of the ideal in
theorem~\ref{algconst}. That is, the residue field of the ideal
constructed also does not contain non-algebraic constants over the
field of constants of the field $K$. We now demonstrate an
application of Theorem~\ref{atomic}.

\begin{theorem}\label{DCFatomic}
Let $K$ be a differential field. Then there exists a differentially
closed field $L$ containing $K$ such that any differentially
finitely generated  $K$-subalgebra of $L$ is locally simple.
\end{theorem}
\begin{proof}

As in the proof of Theorem~\ref{DCFalgconst}, we consider the se
$\{B_\alpha\}$ of all simple differentially finitely generated
algebra over $K$ up to isomorphism. Consider an ideal $\frak p$ of
the ring $R=\otimes_\alpha B_\alpha$, as in Theorem~\ref{atomic}. We
consider $L_1$ to be its residue field. By applying the same
procedure to $L_1$ we obtain $L_2$, for $L_2$ we obtain $L_3$, and
so on. As a result we obtain the ascending chain of fields
$$
K=L_0\subseteq L_1\subseteq\ldots\subseteq L_n\subseteq\ldots
$$
We define the sought-for differential field: $L=\cup_k L_k$. The
differential closedness is proved in the same fashion as in
theorem~\ref{DCFalgconst}, therefore we omit the proof of it.

We claim that the field thus obtained has the indicated property. We
conduct induction on the number of the field. For $L_1$ this is
clear from construction. We show that transition from $k$ to $k+1$.
Suppose that a differentially finitely generated algebra $B$ over
$K$ is contained in $L_{k+1}$; then the algebra $L_k\cdot B$ is
locally simple by the construction of $L_{k+1}$. Let $h$ be an
element $h$ such that $(L_k\cdot B)_h$ is simple. By
Corollary~\ref{splitting_corollary} we obtain that
$$
( L_k\cdot B)_{sh}= L_k\otimes_C (C\cdot B)_{sh},
$$
for some $s$ and $C$. As in Theorem~\ref{atomic}, $C$ can be
considered to be a simple differentially finitely generated
$K$-algebra. Then Proposition~\ref{surjstate} guarantees that
$C\cdot B$ is locally simple, and therefore Proposition~7(b) in
\textsection~10 of Ch.~III in~\cite{Kl} guarantees that the algebra
$B$ is locally simple.
\end{proof}

\subsection{Universal extensions}\label{sec43}

The field constructed in Theorem~\ref{DCFatomic} is already almost
the differential closure. Recall the definition.

\begin{definition}
A differentially closed field $L$ containing $K$ is called a {\it
differential closure} of $K$ if for any differentially closed field
$D$ containing $K$ there exists an embedding of $L$ into $D$ over
$K$.
\end{definition}

Below we construct a certain differential closure in quite a
specific way. Such a construction is a vital necessity on the road
to construction of a splitting field.

Suppose that some well-ordering is defined on the set $\Lambda$
(recall: we consider the ring $R=\otimes_{\alpha\in
\Lambda}B_\alpha$). Then we define a well-ordered family of
differential subrings $\{R_\alpha\}$ as follows:
\begin{enumerate}
\item $R_0=B_0$;

\item $R_{\alpha+1}=R_{\alpha}\otimes B_{\alpha+1}$;

\item $R_\alpha=\cup_{\beta<\alpha} R_\beta$ for limit $\alpha$.
\end{enumerate}
If $\frak p$ is some ideal of the ring $R$, we denote by $\frak
p_\alpha$ its contraction to $R_\alpha$.

\begin{theorem}\label{univers}
Let $\{B_\alpha\}_{\alpha\in\Lambda}$ be some set of simple
differentially finitely generated $K$-algebras, and suppose that
some well-ordering is defined on the set $\Lambda$. Then there
exists a prime differential ideal $\frak p$ of $R=\otimes_\alpha
B_\alpha$ such that any nonzero differential ideal of the ring
$R_{\alpha+1}/\frak p_{\alpha+1}$ contracts to a nonzero ideal of
the ring $R_{\alpha}/\frak p_\alpha$.
\end{theorem}
\begin{proof}
We construct a prime differential ideal $\frak p_\alpha$ of the ring
$R_\alpha$ by transfinite induction. We set $\frak p_0=0$, and at
limit ordinals, $\frak p_\alpha=\cup_{\beta<\alpha}\frak p_\beta$.
It remains to consider the case of non-limit ordinals.

Suppose that $\frak p_\alpha$ has already been constructed; we
construct $\frak p_{\alpha+1}$ in $R_{\alpha+1}$ as a maximal
differential ideal contracting to $\frak p_\alpha$. Then the ideal
$\frak p = \cup_\alpha \frak p_\alpha$ by construction has all the
required properties.
\end{proof}

\begin{proposition}\label{univers_statement}
Let $\frak p$ be the same ideal as in theorem~\ref{univers}. Then
for any differentially closed field $L$ containing $K$ there exists
an embedding of the residue field of the ideal $\frak p$ into $L$
over $K$.
\end{proposition}
\begin{proof}

First we make one simple remark. For any differential subfield $F$
of the differentially closed field $L$ and for any simple
differentially finitely generated algebra $B$ over the field $F$
there exists an embedding of $B$ into $L$ over $F$.

We define the field $K_\alpha$ as the residue field of the ideal
$\frak p_\alpha$. We see that the residue field of $\frak p$ is the
union of the fields $K_\alpha$. We construct an embedding by
transfinite induction. By our remark there exists an embedding of
$R_0$ into $L$, and, consequently, also an embedding of $K_0$ into
$L$. Suppose that $K_\alpha$ is already embedded; then
$$
K_\alpha \cdot (R_{\alpha+1}/\frak p_{\alpha+1})=S^{-1}
R_{\alpha+1}/\frak p_{\alpha+1}
$$
, where $S=R_\alpha\setminus\frak p_\alpha$, is a simple algebra by
the definition of the ideal $\frak p$. Consequently, $K_{\alpha+1}$
can also be embedded. At limit ordinals the embedding is obtained
automatically.
\end{proof}

\begin{corollary}\label{univers_corollary}
Let $\frak p$ be the same ideal as in Theorem~\ref{univers}. Then in
its residue field every differentially finitely generated
$K$-subalgebra is locally simple.
\end{corollary}
\begin{proof}
Let $L$ be the differentially closed field containing $K$ described
in Theorem~\ref{DCFatomic}. Then the residue field of the idea
$\frak p$ can e embedded into $L$. Consequently, the required
property holds.
\end{proof}

We now describe the process of constructing a differential closure.
Subsequently its structure will be very useful for us.

\begin{theorem}\label{DCFunivers}
Let $K$ be a differential field. Then there exist a differentially
closed field $L$ and a well-ordered chain of differential fields
$\{\,L_{\alpha}\,\}$ such that
\begin{enumerate}
\item $L_0=K$;

\item $L=\cup_\alpha L_\alpha$;

\item  $L_{\alpha+1}$ is differentially finitely generated over
$L_\alpha$;

\item every subalgebra of $L_{\alpha+1}$ that is differentially
finitely generated over $L_\alpha$ is locally simple.
\end{enumerate}
Furthermore, the field $L$ is a differential closure of the field
$K$.
\end{theorem}
\begin{proof}

Consider the set $\{\,B_\alpha\,\}$ of all simple differentially
finitely generated algebras over $K$ up to isomorphism. Consider the
ideal $\frak p$ of the ring $R=\otimes_\alpha B_\alpha$, as in
Theorem~\ref{univers}. Let $K_1$ be the residue field of the ideal
$\frak p$. We define $K_{1\,\alpha}$ as the residue field of the
ideal $\frak p_\alpha$. By applying the same procedure to $K_1$, we
obtain $K_2$, for $K_2$ we obtain $K_3$, and so on. As a result we
obtain an ascending chain of differential fields
$$
K= K_0\subseteq K_1\subseteq\ldots\subseteq K_n\subseteq\ldots
$$
We set $ L=\cup_k K_k$. By taking the union of countably many
well-ordered sets we obtain a unified numbering for all the
subfields $K_{k\,\alpha}$. This chain satisfies all the necessary
requirements by Proposition 7(b) in \textsection~10 of
Ch.~\cite{Kl}.

The fact that the field $L$ is differentially closed is proved in
the same fashion as in Theorem~\ref{DCFalgconst}. Since by
construction $K_{k+1}$ is universal over $K_{k}$ in the indicated
sense (Proposition~\ref{univers_statement}), the whole field $L$
also has the same property over $K$.
\end{proof}

\section{Constructed fields}\label{sec5}

This section is devoted to the more detailed analysis of
differential fields constructed in Theorem~\ref{DCFunivers}. We give
the following definition.

\begin{definition}\label{defconstr}
Suppose that we are given some differential field $L$ containing the
field $K$ and a well-ordered chain of differential subfields
$\{L_\alpha\}$ of $L$ such that
\begin{enumerate}
\item $L_0=K$;

\item $L=\cup_\alpha L_\alpha$;

\item  $L_{\alpha+1}$ is differentially finitely generated
over $L_\alpha$;

\item every subalgebra of $L_{\alpha+1}$ that is differentially finitely generated over
$L_\alpha$ is locally simple.
\end{enumerate}
\end{definition}
Then we say that $L$ is constructed over $K$, and call the chain of
differential fields $\{L_\alpha\}$ a {\it construction of the field}
$L$.

Terms is borrowed from \textsection~10.4 in~\cite{Pz}. We begin with
following proposition.

\begin{proposition}\label{prop_constr}

Suppose that we are given some differential field $L$ containing the
field $K$ a well-ordered chain of differential subfields
$\{\L_\alpha\}$ of $L$ such that
\begin{enumerate}
\item $L_0=K$

\item  $L=\cup_\alpha L_\alpha$

\item $L_{\alpha+1}$ is at most countably differentially generated
over $L_\alpha$

\item every subalgebra $L_{\alpha+1}$ that is differentially finitely
generated over $L_\alpha$ is locally simple
\end{enumerate}
Then $L$ is constructed; moreover, the chain $\{L_\alpha\}$ can be
refined to a construction of the field $L$.
\end{proposition}
\begin{proof}

It suffices to learn how to refine to a construction one `storey'
$L_\alpha\subseteq L_{\alpha+1}$. Let $\{x_n\}_{n=0}^\infty$ be a
system of differential generators of $L_{\alpha+1}$ over $L_\alpha$.
We define a chain of differential $L_\alpha$-algebras $B_n$ and
their fields of fractions $F_n$ as follows:
$$
B_n=L_\alpha\{x_1,\ldots,x_n\}.
$$
We claim that the field $F_{n+1}$ over the field $F_n$ has
property~4) in Definition~\ref{defconstr} of a constructed field.
Then, by taking together all such fields for all ordinals, we obtain
a required construction.

Let $B$ be a differentially finitely generated $F_n$-algebra in
$F_{n+1}$. Then it has the form
$$
B=F_n\{z_1,\ldots,z_k\}
$$
for some $z_1,\ldots,z_k\in F_{n+1}$. Thus, $B$ is a localization of
the algebra $B_n\{z_1,\ldots,z_k\}$. The latter is locally simple by
condition~4) in the statement of the proposition. Since a
localization of a locally simple algebra is locally simple, we
obtain what is required.
\end{proof}

We now state conditions that guarantee that a field is constructed
over a subfield. We introduce several auxiliary definitions, also
borrowed from~\textsection~10.4 in~\cite{Pz}.

\begin{definition}
Let $L$ be a constructed differential field over $K$. Then by
conditions~3) and~4) in Definition~\ref{defconstr} for any ordinal
$\alpha$ there exists a simple differentially finitely generated
$K$-algebra $B_\alpha\subseteq L_{\alpha+1}$ such that the field
$L_{\alpha+1}$ is the field of fractions of the algebra
$L_\alpha\cdot B_\alpha$. The family of differential $K$-algebras
$\{B_\alpha\}$ is called a {\it generating set} for the constructed
field $L$.
\end{definition}

We should point out that, first, a generating family can be chosen
in different ways. Second, every ring in a generating family is a
differentially finitely generated algebra over the ground field $K$.
Note that the field $L_\alpha$ can be reconstructed from any
generating family as follows:
$$
L_{\alpha+1}=L_\alpha\langle B_\alpha \rangle,
$$
and for limit ordinals
$$
L_\alpha=K\langle\{B_\beta\}_{\beta<\alpha} \rangle.
$$

\begin{proposition}\label{prop_package}
Let $L$ be a constructed field over $K$, and $\{\,B_\alpha\,\}$ some
generating family. Then for every $B_\alpha$ there exist a finite
set $\{\,B_{\alpha_1},\ldots,B_{\alpha_t}\,\}$ with the condition
$\alpha_i<\alpha$ and an element  $s\in F_\alpha \cdot B_\alpha$
such that
$$
(L_\alpha\cdot B_\alpha)_s=L_\alpha\otimes_{F_\alpha}(F_\alpha\cdot
B_\alpha)_s,
$$
where
$$
F_\alpha=K\langle B_{\alpha_1},\ldots,B_{\alpha_t}\rangle
$$
and the algebra $(F_\alpha\cdot B_\alpha)_s$ is simple.
\end{proposition}
\begin{proof}

The algebra $L_\alpha\cdot B_\alpha$ is locally simple;
consequently, there exists an element $h\in L_\alpha\cdot B_\alpha$
such that the algebra
$$
(L_\alpha\cdot B_\alpha)_h
$$
is simple. Then from the Corollary~\ref{splitting_corollary} of the
``splitting lemma'' there exist a differentially finitely generated
$K$-algebra $C\subseteq L_\alpha\cdot B_\alpha$ and an element
$s'\in C\cdot B_\alpha$ such that
$$
(L_\alpha\cdot B_\alpha)_{s'h}= L_\alpha \otimes_{C} (C\cdot
B_\alpha)_{s'h}.
$$
We denote the product $s'h$ by $s$. Since the ring $C$ is
differentially finitely generated over $K$, there exists a finite
set of $K$-algebras
$$
\{B_{\alpha_1},\ldots,B_{\alpha_t}\}
$$
such that $C$ belongs to the field
$$
F_\alpha=K\langle B_{\alpha_1},\ldots,B_{\alpha_t}\rangle,
$$
and we have the decomposition
$$
(L_\alpha\cdot B_\alpha)_s= L_\alpha \otimes_{F_\alpha}
(F_\alpha\cdot B_\alpha)_s.
$$
Since the right-hand side of the equation is a simple differential
ring, it follows from Proposition~\ref{surjstate} that
$(F_\alpha\cdot B_\alpha)_s$ is also simple.
\end{proof}

Suppose that some generating family of $K$-algebras
$\{\,B_\alpha\,\}$ is fixed for a differential field $L$.We define
the notion of a package of a ring $B_\alpha$.

\begin{definition}\label{defpackage}
Let $L$ be a constructed differential field over $L$ and
$\{B_\alpha\}$ some generating family. We give the definition by
induction. The {\it package $\pi_0$ of the algebra $B_0$} is
declared to be the empty set. Suppose that the package is defined
for $B_\beta$ for all $\beta<\alpha$; we define the package of
$B_\alpha$. Let $\{B_{\alpha_1},\ldots,B_{\alpha_t}\}$ be the family
of $K$-algebras for the algebra $B_\alpha$ as in
Proposition~\ref{prop_package}. Then the {\it package $\pi_\alpha$}
is defined to be the union
$$
\{B_{\alpha_1},\ldots,B_{\alpha_t}\}\cup\bigcup_i \pi_{\alpha_i}.
$$
\end{definition}

By induction we obtain that the package of any $K$-algebra in a
generating family consists of finitely many $K$-algebras.

\begin{definition}
The {\it package field of the algebra} $B_\alpha$ is defined to be
the differential field generated by all the algebras in the package:
$$
F_\alpha=K\angle\pi_\alpha \rangle
$$
\end{definition}

Note that packages are not uniquely defined. The definition of a
package and Proposition~\ref{prop_package} imply the following
lemma.

\begin{lemma}[``on package splitting'']\label{lemma_package}
Suppose that $L$ is a constructed differential field over $K$,
$\{\,B_\alpha\,\}$ is some generating family, and for each $\alpha$
some package $\pi_\alpha$ is fixed. Then for any $\alpha$ there
exists an element $s\in F_\alpha \cdot B_\alpha$ such that
$$
(L_\alpha\cdot B_\alpha)_s=L_\alpha\otimes_{F_\alpha}(F_\alpha\cdot
B_\alpha)_s.
$$
Furthermore, the algebra $(F_\alpha\cdot B_\alpha)_s$  is simple.
\end{lemma}

\begin{proposition}\label{pocketcontent}
Let $L$ be a constructed differential field over $K$ with some
family of generating $K$-algebras $\{\,B_\alpha\,\}$. Suppose that
for every $\alpha$ some package $\pi_\alpha$ is fixed and for some
$\alpha$ a subfield $F$ contains the field $F_\alpha$, and let $s$
be the same elements as in the lemma on package splitting. Then the
algebra $(F\cdot B_\alpha)_s$ is simple.
\end{proposition}
\begin{proof}
By Lemma~\ref{lemma_package}, for this element $s$ we have
$$
(L_\alpha\cdot B_\alpha)_s=L_\alpha\otimes_{F_\alpha}(F_\alpha\cdot
B_\alpha)_s.
$$
Since $F$ contains $F_\alpha$, we have
$$
(L_\alpha\cdot
B_\alpha)_s=L_\alpha\otimes_{F}F\otimes_{F_\alpha}(F_\alpha\cdot
B_\alpha)_s=L_\alpha\otimes_{F}(F\cdot B_\alpha)_s.
$$
It follows from the choice of the element $s$ that the algebra on
the left-hand side of the last equation is simple; therefore by
Proposition~\ref{surjstate} the algebra $(F\cdot B_\alpha)_s$ is
also simple.
\end{proof}

\begin{definition}

Let $L$ be a constructed differential field over $K$, and
$\{B_\alpha\}$ some generating family. Suppose that for every
$K$-algebra $B_\alpha$ in the generating family
$\{\,B_\alpha\,\}_{\alpha\in \Lambda}$ some package $\pi_\alpha$ is
fixed. Consider a subfamily of $K$-algebras
$\{\,B_\beta\,\}_{\beta\in X}$ ($X\subseteq\Lambda$) such that with
any ring $B_\beta$, $X\subseteq\Lambda$, such that, with each
$K$-algebra $B_\beta$, this subfamily contains its entire package.
Such a family is called {\it package closed} or, more simply, {\it
closed}.
\end{definition}

This term is used in \textsection~10.4~\cite{Pz}. Recall that the
differential finite generatedness of algebras in a general family is
a part of their definition.

\begin{definition}
Suppose that the conditions of the preceding definition hold, and
suppose that a differential subfield $F\subseteq L$ is generated by
some closed family of algebras $\{B_\beta\}$. Then such a field is
also called {\it closed}.
\end{definition}

A subset $X$ is a well-ordered set with respect to the ordering
induced from the set $\Lambda$. Then we can define a family of
subfields $\{\,L'_\beta\,\}$ as follows:
$$
L'_\beta=K\langle\{\,B_{\theta}\,\}_{\substack{\theta<\beta\\\theta\in
X}}\rangle.
$$
The importance of closed subfields is emphasized by the following
two propositions (see their logical originals in~\cite{Pz}, Ch.~10,
Propositions~10.15, 10.17.

\begin{proposition}\label{prop_closed_constr}
Let $L$ be a constructed differential field over $K$ with some
family of generating algebras $\{\,B_\alpha\,\}$. Suppose that for
each $\alpha$ some package is fixed, and suppose that some
differential subfield $F$ of the field $L$ is closed. Then the field
$F$ is constructed over $K$, and the chain of fields
$\{\,L'_\beta\,\}$ is a construction.
\end{proposition}
\begin{proof}
We need to show that property~4) in definition~\ref{defconstr} holds
for the pair of fields $L'_\beta\subseteq L'_{\beta+1}$. By
definition, $L'_{\beta+1}=L'_\beta\langle B_\gamma\rangle$ for some
$\gamma\in X$, for every $\theta\in X$ strictly smaller than
$\gamma$ we have $B_\theta\subseteq L'_\beta$. By Proposition~7(b)
in \textsection~10 of Ch.~III in~\cite{Kl} it suffices to show that
the algebra $L'_\beta\cdot B_\gamma$ is locally simple. However, the
field $F$ is closed and therefore contains $F_\gamma$. By definition
of a package, $F_\gamma$ is contained in $L'_\beta$, since it is
generated by tings of the form $B_\theta$, where $\theta\in X$ and
strictly smaller than $\gamma$. Then Proposition~\ref{pocketcontent}
implies what is required.
\end{proof}

\begin{proposition}\label{pocketclose}
Let $L$ be a constructed differential field over $K$ with some
family of generating algebras $\{\,B_\alpha\,\}$. Suppose that for
each $\alpha$ some package $\pi_\alpha$ is fixed and some
differential  subfield $F$ of the field $L$ is closed. Then the
field $L$ is constructed over $F$, and the chain of fields
$\{\,F{<}L_\alpha{>}\,\}$ is a construction.
\end{proposition}
\begin{proof}

We need to show that property~4) in definition~\ref{defconstr} holds
for the pair of fields $F\langle L_\alpha\rangle\subseteq F\langle
L_{\alpha+1}\rangle$. By Proposition~7(b) in \textsection~10 of
Ch.~III in~\cite{Kl} it suffices to show that the algebra
$$
F\langle L_\alpha\rangle\cdot B_\alpha
$$
is locally simple. by definition the algebra $L_\alpha\cdot
B_\alpha$ is locally simple. Let $s$ be an element such that algebra
$$
(L_\alpha\cdot B_\alpha)_s
$$
is simple. We shall prove that algebra
$$
(F\langle L_\alpha\rangle\cdot B_\alpha)_s
$$
is also simple. Since $F$ is constructed by
Proposition~\ref{prop_closed_constr}, we use induction on $\beta$ to
show that the algebra
$$
(L'_\beta\langle L_\alpha\rangle \cdot B_\alpha)_s
$$
is simple. Since a limit of simple differential ring is a simple
differential ring, it sufficient to consider non-limit
ordinals.

To begin with, we observe that the equation
$$
L'_\beta\langle L_\alpha\rangle=\Qt(L'_\beta\cdot B_\alpha)
$$
holds for any $\beta$. Let $\widetilde{F}=\Qt(L'_\beta\cdot
B_\alpha)$. Since the field $\Qt(\widetilde{F}\cdot B_\alpha)$
contains the package of $B_{\beta+1}$, by
Proposition~\ref{pocketcontent} for some $h\in\widetilde{F}\cdot
B_{\beta+1}$ the algebra
$$
(\Qt(\widetilde{F}\cdot B_\alpha)\cdot B_{\beta+1})_h
$$
is simple. We set $S=(\widetilde{F}\cdot B_\alpha)_s\setminus
\{\,0\,\}$; then
\begin{multline*}
(\Qt(\widetilde{F}\cdot B_\alpha)\cdot
B_{\beta+1})_h=(S^{-1}((\widetilde{F}\cdot B_\alpha)_s)\cdot
B_{\beta+1})_h=\\
=S^{-1}(((\widetilde{F}\cdot B_\alpha)_s\cdot
B_{\beta+1})_h)=S^{-1}(((\widetilde{F}\cdot B_{\beta+1})_h\cdot
B_\alpha)_s)
\end{multline*}
Therefore the algebra
$$
S^{-1}(((\widetilde{F}\cdot B_{\beta+1})_h\cdot B_\alpha)_s)
$$
is simple.  But since $(\widetilde{F}\cdot B_\alpha)_s$ is simple by
induction, the set $S$ cannot intersect any differential ideal in
any differential ring. Therefore the ring
$$
((\widetilde{F}\cdot B_{\beta+1})_h\cdot B_\alpha)_s
$$
is also simple, and together with it, so is also
$$
(L'_{\beta+1}\langle L_\alpha\rangle\cdot B_\alpha)_s
$$
as its localization.
\end{proof}

\begin{proposition}\label{constatomic}
Let $L$ be a constructed differential field over $K$. Then any
differentially finitely generated $K$-subalgebra of $L$ is locally
simple.
\end{proposition}
\begin{proof}
Any constructed differential field $L$ has the following property:
for any differential homomorphism of the field $K$ into a
differentially closed field there exists an extension of it to $L$.
Therefore our field can be embedded into the field in
Theorem~\ref{DCFatomic}, whence the desired result follows.
\end{proof}

\begin{proposition}\label{constrclos}
Let $L$ be a constructed differential field over $K$ with some
family of generating $K$-algebras $\{\,B_\alpha\,\}$. Suppose that
for each $\alpha$ some package $\pi_\alpha$ is fixed, and let $F$ be
some subfield of $L$ that is differentially finitely generated over
$K$. Then there exists a closed subfield $F'$ containing $F$ that is
differentially finitely generated over $K$.
\end{proposition}
\begin{proof}
Since the field $F$ is differentially finitely generated, it is
contained in the field generated by some rings
$B_{\alpha_1},\ldots,B_{\alpha_k}$. We consider as the differential
field $F'$ the field generated by these algebras and their packages.
We now see that it satisfies the required properties.
\end{proof}

\begin{proposition}\label{FiniteConstr}
Let $L$ be a constructed differential field over $K$. Then $L$ is
constructed over any subfield of it that is differentially finitely
generated over $K$.
\end{proposition}
\begin{proof}
Let $F$ be differentially finitely generated subfield of $L$. Then
by Proposition~\ref{constrclos} there exists a closed subfield $F'$
containing $F$ that is differentially finitely generated over $K$.
Consequently, by Proposition~\ref{pocketclose} $L$ is constructed
over $F'$. In turn, it follows from Proposition~7(b)
in~\textsection~10 of~Ch.~III in~\cite{Kl} that $F'$ is constructed
over $F$, and therefore $L$ is also constructed over $F$.
\end{proof}

The following result in model theory is known as Ressayre's
theorem(see~\cite[chapter~10, sec.~4, theor.~10.18]{Pz}); we present
its algebraic analogue.

\begin{theorem}\label{ConstructiveUnique}
Let $L$ and $F$ be differentially closed fields that are constructed
over $K$. Then they are isomorphic.
\end{theorem}
\begin{proof}
Let construction of the field $L$ and $F$ be $\{\,L_\alpha\,\}$ and
$\{\,F_\beta\,\}$, respectively. We also assume that some families
of generating $K$-algebra are fixed on them and their packages are
defined. Then we consider the following set:
$$
\Sigma=\{\, (L',F',f')  \mid L'\subseteq L,\: F'\subseteq F,\:
f\colon L'\to F' \,\},
$$
where $L'$ and $F'$ are closed subfields and $f'$ is differential
isomorphism between them. On this set we introduce the following
partial order relation:
$$
(L',F',f') \leqslant (L'',F'',f'')\: \Leftrightarrow \: L'\subseteq
L'',\: F'\subseteq F'',\: f''\mid_{L'}=f'
$$
The set $\Sigma$ is not empty, since it contains the element
$(K,K,Id)$, and satisfies the hypothesis of Zorn's lemma.
Consequently, there exists a maximal element
$(\widehat{L},\widehat{F},\widehat{f})$. We claim that
$\widehat{L}=L$ and $\widehat{F}=F$. Suppose the opposite, for
example, $\widehat{L}\neq L$.

Proposition~\ref{pocketclose} guarantees that the field $L$ is
constructed over $\widehat{L}$ (respectively, $F$ over
$\widetilde{F}$), and a generating set of rings over $\widetilde{L}$
is the same as over $K$, and therefore the packages are the same. By
the assumption there exists a simple ring $B_\alpha\subseteq L$ that
is not contained in $\widehat{L}$. by the definition of
constructedness there exists an element $s\in \widehat{L}\cdot
B_\alpha$ such that $(\widehat{L}\cdot B_\alpha)_s$ is a simple
algebra. Then by the differential closedness of $F$ the differential
homomorphism $\widehat{f}$ can be extended to the field of fraction
of the algebra $\widehat{L}\cdot B_\alpha$. We denote this field by
$\widehat{L}_1$. It follows from Proposition~\ref{constrclos} that
there exists a closed differentially finitely generated field
$\widehat{F}_1$ containing the image of $\widehat{L}_1$. In turn,
from Proposition ~\ref{constatomic} and Proposition~7(b)
in~\textsection~10 of~Ch.~III in~\cite{Kl} we obtain that
$\widehat{F}_1$ is the field of fractions of a simple differentially
finitely generated algebra over the image of $\widehat{L}_1$. Then
by the differential closedness of $L$ there exists a reverse
embedding of the field $\widehat{F}_1$ into $L$. Again by
Proposition~\ref{constrclos} there exists a closed differentially
finitely generated field $\widehat{L}_2$ containing the image of
$\widehat{F}_1$. Then we apply this construction to the field
$\widehat{L}_2$, and so on. The scheme of extending the isomorphism
$\widehat{f}$ is presented in the following diagram:
$$
\xymatrix{
    \widehat{L}\ar[d]^{\widehat{f}}&\widehat{L}_1\ar[d]&\widehat{L}_2\ar[d]&\ldots&\widehat{L}_n\ar[d]&\ldots\\
    \widehat{F}&\widehat{F}_1\ar[ur]&\widehat{F}_2\ar[ur]&\ldots&\widehat{F}_n\ar[ur]&\ldots
}
$$
As a result we obtain two chains of subfields: $\widehat{L}_k$ in
$L$, and $\widehat{F}_k$ in $F$, and each element of these chains is
closed. Then we define the fields as $L'=\cup_k \widehat{L}_k$ and
$F'=\cup_k\widehat{F}_k$ and we see that they are closed and by
construction the differential homomorphism $\widehat{f}$ can be
extended to them. We have obtained a contradiction with the
maximality.
\end{proof}

\section{Uniqueness theorem}\label{sec6}

Recall that all differential rings are assumed to be algebras over
some differential field $K$. Theorem~\ref{DCFuniquetheorem} below
contains a technique similar to that presented in Proposition~18.1
in~Ch.~18 of~\cite{Pz}. It is on this technique that the
entire~\textsection~7 is based. To begin with, we need a certain
variant of the splitting lemma. Apparently, the following lemma is
directly related to Corollary~16.7 and Theorem~16.8 in~Ch.~16
of~\cite{Pz}; possibly, it is their algebraic analogue.

\begin{lemma}\label{lemma_open}
Let $A$ and $B$ be arbitrary differential $K$-algebras without
nilpotents and suppose that there exists an element $h\in A\otimes
B$ such that $(A\otimes B)_h$ is a simple differential algebra. Then
$B$ is locally simple.
\end{lemma}
\begin{proof}

The element $h$ has the form $\sum a_i\otimes b_i$. We can assume
that the elements $a_i$ are linearly independent over $K$. Since
$A\otimes B$ is locally simple, for any nonzero differential ideal
$\frak a$ we have
$$
((A\otimes B)/\frak a)_h=0
$$
In particular, for any nonzero prime differential ideal $\frak p$ of
$B$,
$$
(A\otimes (B/\frak p))_h=0.
$$
Since  $A$ and $B/\frak p$ are algebras over the field $K$ (of
characteristic zero) that do not contain nilpotents, it follows that
$A\otimes (B/\frak p)$ also does not contain nilpotents (see~
\cite[lemma~A.16]{SW}). The condition of the ring being equal to
zero means that $h^n=0$ in $A\otimes (B/\frak p)$, and therefore
also $h=0$, that is, the elements $b_i$ belong to $\frak p$.
Consequently, the differential spectrum of the ring $B_{b_i}$
consists of a single element for every $i$. In view of the absence
of nilpotents this means what is required.
\end{proof}

Let $L$ be some differential closure of $K$, and let field $D$ be a
constructed differential closure (which exists by
Theorem~\ref{DCFunivers}). Then we see that the field $L$ can be
embedded into $D$. Our next task is to find a construction of the
differential closure contained in the constructed one.

\begin{theorem}\label{DCFuniquetheorem}
Let $L$ be a differential closure of the field $K$. Then $L$ is
constructed.
\end{theorem}
\begin{proof}

The arguments before the statement of the theorem show that we can
consider the field $L$ to be embedded into some constructed
differential closure $D$ with construction $D_\alpha$. In order to
construct a construction of $L$, it seems natural to choose as such
the corresponding intersections $L_\alpha=L\cap D_\alpha$. However,
it may happen that the fields $L_\alpha$ do not form a construction.
In order to obtain a correct result, we need to touch up slightly
the construction of $D$. Namely, we need to construct a system of
subfields $d'_\alpha$ with the following conditions:
\begin{enumerate}
\item $D_\alpha\subseteq D'_\alpha$;

\item $D'_\alpha$ is closed;

\item $D'_{\alpha+1}$ is at most countable differentially generated
over $D'_{\alpha}$;

\item $L\cdot D'_\alpha=L\otimes_{L_\alpha}D'_\alpha$, where $L_\alpha=L\cap
D'_\alpha$;

\item $L_{\alpha+1}$ is at most countable differentially generated
over $L_\alpha$.
\end{enumerate}

These conditions are preserved under passing to inductive limits;
therefore it suffices to construct such fields for non-limit
ordinals. We assume that some generating set of $K$-algebras
$\{B_\alpha\}$ is chosen in $D$ and packages are defined for them.
Thus, the notion of being closed is defined.

We construct a chain of differential fields $\overline{D}_k$ as
follows. We set $\overline{D}_1=D'_\alpha\langle D_\alpha\rangle$,
it is differentially finitely generated over $D'_\alpha$.  By
construction, $\overline{D}_1$ is closed; however, it is not
guaranteed that $\overline{D}_1$ and $L$ are linearly disjoint, that
is, generally speaking,
$$
L\cdot \overline{D}_1\neq L\otimes_{L\cap \overline{D}_1}
\overline{D}_1.
$$
By Corollary~\ref{splitting_field_corollary} (using the field $L\cap
D'_\alpha$ instead of $K$) we obtain that there exist a field $C_1$
differentially finitely generated over $L\cap D'_\alpha$ and an
element $s\in C_1\cdot \overline{D}_1$ such that
$$
(L\cdot \overline{D}_1)_s=L\otimes_{C_1}(C_1\cdot \overline{D}_1)_s.
$$
In this case, by setting $\overline{D}_2=\Qt(C_1\cdot
\overline{D}_1)$, we have
$$
L\cdot \overline{D}_2=L\otimes_{C_1}\overline{D}_2.
$$
Then by the definition of the tensor product we have $C_1=L\cap
\overline{D}_2$ and this field is also differentially finitely
generated over $L\cap D'_\alpha$. But now the field $\overline{D}_2$
is not necessarily closed. Proposition~\ref{constrclos} guarantees
the existence of a closed field $\overline{D}_3$ differentially
finitely generated over $D'_\alpha$ and containing $\overline{D}_2$.
By repeating for $\overline{D}_3$ the procedure carried out for
$\overline{D}_1$ we obtain a field $\overline{D}_4$ containing
$\overline{D}_3$ such that
$$
L\cdot \overline{D}_4=L\otimes_{C_2}\overline{D}_4,
$$
and $C_2=L\cap \overline{D}_4$ is differentially finitely generated
over $L\cap D'_\alpha$. Consequently, $C_1\subseteq C_2$. By
proceeding in similar fashion, we construct a chain of differential
fields $\overline{D}_k$ and $C_k$ with the following conditions:
\begin{enumerate}
\item $C_k$ is differentially finitely generated over $L\cap
D'_\alpha$;
\item $\overline{D}_k$ is differentially finitely generated over
$D'_\alpha$;
\item $C_{n}=L\cap \overline{D}_{2n}$;
\item $L\cdot \overline{D}_{2n}=L\otimes_{C_n}\overline{D}_{2n}$;
\item $\overline{D}_{2n+1}$ is closed.
\end{enumerate}
We set $D'_{\alpha+1}=\cup_k \overline{D}_k$. Then by property~5)
this field is closed. By passing to the direct limit (see details
in~\cite[Ch.~2, Exercise~20]{AM}) we obtain
$$
L\cdot D'_{\alpha+1}=L\otimes_{L_{\alpha+1}} D'_{\alpha+1},
$$
where $L_{\alpha+1}=L\cap D'_{\alpha+1}$, and $L_{\alpha+1}=\cup_k
C_k$. As we see, all five properties are satisfied.

We claim that the chain of subfields $\{\,L_\alpha\,\}$ satisfies
the properties of Proposition~\ref{prop_constr}. For that we only
need to verify property~4). Let $B$ be an arbitrary differentially
finitely generated algebra over $L_\alpha$; then algebra
$D'_\alpha\cdot B$ is differentially finitely generated over
$D'_\alpha$. By the closedness of $D'_\alpha$ it follows from
Propositions~\ref{pocketclose} and~\ref{constatomic} that
$D'_\alpha\cdot B$ is locally simple. We now consider the following
chain of inclusions:
$$
D'_\alpha\otimes_{L_\alpha}B\subseteq
D'_\alpha\otimes_{L_\alpha}L=D'_\alpha\cdot L
$$
Consequently,
$$
D'_\alpha\cdot B=D'_\alpha\otimes_{L_\alpha}B.
$$
Since $D'_\alpha\otimes_{L_\alpha}B$ is locally simple, it follows
from Lemma~\ref{lemma_open} (for $L_\alpha$ instead of $K$) that $B$
is also locally simple.
\end{proof}

\begin{corollary}\label{DCFuniversUnique}
A differential closure is unique up to isomorphism, and it is
constructed.
\end{corollary}

\section{Splitting fields}\label{sec7}

\subsection{Abstract splitting fields}\label{sec71}

As above, all differential rings are assumed to be algebras over
some differential field $K$. We consider an arbitrary family of
differentially finitely generated $K$-algebras $\{\,B_\alpha\,\}$
(not necessarily simple ones), and let $L$ be some differential
field. We consider all possible differential homomorphisms
$B_\alpha\to L$ for all possible $\alpha$. Thus some family of
differential subrings of $L$ is defined. If $L$ is the smallest
differential field containing this family, then we shall say that
$L$ {{\it is generated by the set of rings} $\{\,B_\alpha\,\}$. For
the convenience of subsequent notation we denote by
$\overline{B_\alpha}$ the image of $B_\alpha$ in $L$ under some
differential homomorphism.

\begin{definition}\label{def_split}
We say that $L$ is a {\it splitting field} of a family of
differentially finitely generated algebras $\{B_\alpha\}$ over $K$
if the following conditions hold:
\begin{enumerate}
\item $L$ is sufficiently large: for any $\alpha$ and for any
locally maximal differential ideal $\frak m$ of $L\otimes B_\alpha$
we have $(L\otimes B_\alpha)/\frak m=L$ (this means as isomorphism
onto the first factor of the tensor product);
\item $L$ is not too large: the field $L$ is generated by the family
$\{B_\alpha\}$;
\item $L$ is universal: for any differential field $L'$ satisfying
properties~1) and~2) there exists an embedding of $L$ into $L'$ over
$K$.
\end{enumerate}
\end{definition}

Note that the differential closure becomes a splitting field for all
differentially finitely generated algebras over the field $K$.

\begin{proposition}\label{SplitDescriptor}
Let $\{\,B_\alpha\,\}$ be some family of differentially finitely
generated $K$-algebra, and suppose that a differential field $L$ has
the following properties:
\begin{enumerate}
\item for any $\alpha$ and for any locally maximal differential
ideal $\frak m$ of $L\otimes B_\alpha$ we have $(L\otimes
B_\alpha)/\frak m=L$;
\item the field $L$ is generated by the family
$\{\,B_\alpha\,\}$;
\item $L$ is constructed over $K$.
\end{enumerate}
Then $L$ is a splitting field of the family  $\{\,B_\alpha\,\}$ over
$K$.
\end{proposition}
\begin{proof}

We need to show that constructedness implies property~3) of the
definition of a splitting field. Indeed, let $L'$ be an arbitrary
field with properties~1) and~2); then we denote by $\overline{L}'$
its differential closure. The constructedness implies that $L$
embeds into $\overline{L}'$. It remains to show that the image of
$L$ is contained in $L'$. Since $L$ is generated by the rings
$\overline{B}_\alpha$, it is sufficient to show that all these rings
are contained in $L'$. The algebra $L'\cdot \overline{B}_\alpha$ is
locally simple by Proposition~\ref{constatomic}; but then by the
definition of $L'$ this algebra coincides with $L'$, as required.
\end{proof}

We now show that for any family of differentially finitely generated
$K$-algebras splitting fields exist and have very special form.

\begin{proposition}\label{ExistSplitClose}
For any set of differentially finitely generated $K$-algebras
$\{\,B_\alpha\,\}$ there exists a constructed splitting field, and
its differential closure is the differential closure of the field
$K$.
\end{proposition}
\begin{proof}

Consider a family of simple differentially finitely generated
$K$-algebras of the following form:
$$
F_1=\{\,(B_\alpha/\frak m)_s\,\},
$$
where $\frak m$ is an arbitrary locally maximal differential ideals
of $B_\alpha$ and $s$ is the corresponding element such that
$(B_\alpha/\frak m)_s$ is simple. We well-order the set $F_1$. Let
this family be well-ordered. Then in the ring
$$
R_1=\otimes_{B\in F_1} B
$$
we choose an ideal $\frak p_1$ as in Theorem~\ref{univers}.The
residue field of the ideal $\frak p_1$ is denoted by $L_1$. Then we
need to repeat the procedure for the family of algebras
$\{\,L_1\otimes B_\alpha\,\}$ instead of $\{\,B_\alpha\,\}$. We
obtain a field $L_2$ containing $L_1$, and so on. We set $L=\cup_k
L_k$. We claim that the field $L$ is the required one. We observe at
once that it is constructed.

We need to show conditions~1) and~3) of Definition~\ref{def_split}
hold. Let us show property~1). Let
$$
(L\cdot \overline{B}_\alpha)_s=(L\otimes B_\alpha)_s/\frak m
$$
be some simple differential algebra. Then by
Corollary~\ref{splitting_corollary} there exist a differentially
finitely generated $K$-algebra $C\subseteq L$ and an element $h\in
C\cdot \overline{B}_\alpha$ such that
$$
(L\cdot \overline{B}_\alpha)_{sh}=L\otimes_{C}(C\cdot
\overline{B}_\alpha)_{sh}.
$$
Since $C$ is differentially finitely generated, we have $C\subseteq
L_k$ for some $k$. Therefore also,
$$
(L\cdot \overline{B}_\alpha)_{sh}=L\otimes_{L_k}(L_k\cdot
\overline{B}_\alpha)_{sh}.
$$
But the algebra $(L_k\cdot \overline{B}_\alpha)_{sh}$ embeds into
$L_{k+1}$ by construction, and therefore also into $L$. That is, the
algebra $(L\cdot \overline{B}_\alpha)_{sh}$ embeds into $L$ and
therefore is also isomorphic to it, as required.

Since the field $L$ is constructed over $K$, the differential
closure of $L$ is the differential closure of $K$. It now follows
from Proposition~\ref{SplitDescriptor} that $L$ is precisely the
sought-for field.
\end{proof}

\begin{proposition}\label{ExtAutoSplit}
For any set of differentially finitely generated $K$-algebras
$\{\,B_\alpha\,\}$, a splitting field is unique, and any
differential automorphism of the field $K$ can be extended to a
differential automorphism of the splitting field.
\end{proposition}
\begin{proof}

Let $L$ and $L'$ be two splitting fields. Then it follows from
Proposition~\ref{ExistSplitClose} and
Corollary~\ref{DCFuniversUnique} that the fields
$$
\overline{L}=\overline{L}'=\overline{K},
$$
can be identified by the following maps:
$$
\xymatrix@R=4pt@C=10pt{
    {\overline{L}}\ar[r]&{\overline{K}}&{\overline{L}'}\ar[l]\\
    {L}\ar@{}[u]|{\cup}&{}&{L'}\ar@{}[u]|{\cup}\\
    {K}\ar@{=}[r]\ar@{}[u]|{\cup}&{K}\ar@{}[uu]|{\bigcup}&{K}\ar@{=}[l]\ar@{}[u]|{\cup}\\
}
$$
It remains only to verify that the fields $L$ and $L'$ coincide
under this identification. Let us prove the inclusion $L\subseteq
L'$. Since the family of differential rings of the form
$\overline{B}_\alpha$ generates $L$ over $K$, it suffices to show
that any algebra $B$ of the form indicated above is contained in
$L'$. Indeed, consider the algebra $L'\cdot B$; then it follows from
equality $\overline{L}=\overline{K}$ and
Proposition~\ref{constatomic} that this algebra is locally simple
and, consequently, $B$ is contained in $L'$ by definition. The
reverse inclusion is verified in similar fashion.

Note that any differential automorphism of the field $K$ can be
extended to a differential automorphism of $\overline{K}$; this is
precisely the theorem of uniqueness of the differential closure. But
then, by identifying the fields by the isomorphism constructed, we
reduce the problem to verifying the uniqueness of a splitting field.
But since this has already been proved, the proposition is also
proved.
\end{proof}

It follows from the definition of a splitting field of a family
$\{\,B_\alpha\,\}$ that for any locally maximal differential ideal
$\frak m$ of $B_\alpha$ the ring $B_\alpha/\frak m$  can be embedded
into the splitting field. We now show that a somewhat more general
result is true.

\begin{proposition}\label{prop_tensor}
Let $L$ be the splitting field over $K$ of a family of
differentially finitely generated $K$-algebras $\{\,B_\alpha\,\}$.
Then for any locally maximal differential ideal $\frak m$ of the
ring
$$
L\otimes B_{\alpha_1}\otimes\cdots\otimes B_{\alpha_n}
$$
its residue field coincides with $L$. As a corollary, any locally
simple differential algebra of the form
$$
(B_{\alpha_1}\otimes\cdots\otimes B_{\alpha_n})/\frak m
$$
embeds into $L$.
\end{proposition}
\begin{proof}
We use induction on the number of factors. For $n=1$ the condition
that is being proved is contained in the definition of a splitting
field. We now suppose that for $n$ everything has already been
proved. We consider the equality
$$
(L\otimes B_{\alpha_1}\otimes\cdots\otimes B_{\alpha_{n+1}})/\frak
m=((L\otimes B_{\alpha_1}\otimes\cdots\otimes B_{\alpha_{n}})/\frak
m^c\otimes B_{\alpha_{n+1}})/\frak m
$$
It follows from Proposition~7(b) in~\textsection~10 of~Ch.~III
in~\cite{Kl} that $\frak m^c$ is a locally maximal differential
ideal and therefore by induction the latter ring is equal to
$$
(L\otimes B_{\alpha_{n+1}})/\frak m=L,
$$
as required.

We denote an algebra of the form $(B_{\alpha_1}\otimes\ldots\otimes
B_{\alpha_n})/\frak m$ by $B$; then for some $s$ in $B$ the algebra
$B_s$ is simple. Therefore for any maximal differential ideals
$\frak m$ we have the equality
$$
(L\otimes B_s)/\frak m=L,
$$
which ensures an embedding of $B$ into $L$.
\end{proof}

\begin{proposition}\label{SplitScalExt}
Let $F$ be a differential subfield of the splitting field $L$ of a
family of differentially finitely generated $K$-algebras
$\{\,B_\alpha\,\}$ over $K$ such that $L$ is constructed over $F$.
Then $L$ is the splitting field over $F$ of the family $\{\,F\otimes
B_\alpha\,\}$.
\end{proposition}
\begin{proof}

The first part of the definition follows from the equality
$$
L\otimes_{F}F\otimes B_\alpha=L\otimes B_\alpha.
$$
Since $L$ is generated by the images of the $B_\alpha$, it follows
that, a fortiori, $L$ is generated by the images of the  $F\otimes
B_\alpha$. Then property~3) of Definition~\ref{def_split} follows
from Proposition~\ref{SplitDescriptor}.
\end{proof}

We denote the group of differential automorphisms of $L$ over $K$ by
$\diffaut(L/K)$. For an arbitrary group $H$ of differential
automorphisms of $L$ over $K$, we denote by $L^H$ the fixed-element
subfield with respect to $H$.

\begin{proposition}\label{InvariantSplit}
For any splitting field $L$ over $K$ we have the equation
$$
K=L^{\diffaut (L/K)}.
$$
\end{proposition}
\begin{proof}

Let $x\in L\setminus K$; then this element belongs to some
subalgebra of the form
$$
B=\overline{B}_{\alpha_1}\cdot \ldots\cdot \overline{B}_{\alpha_n}.
$$
Since this algebra is differentially finitely generated and $L$ is
constructed, Proposition~\ref{constatomic} guarantees that $B$ is
locally simple, that is, for some  $s\in B$ the algebra $D=B_s$ is
simple. Since $D$ is a Ritt algebra, the ring $D\otimes D$  has no
nilpotents (see~\cite{SW}, Lemma A.16). Let $\frak m$ be a maximal
differential ideal of the ring
$$
(D\otimes D)_{1\otimes x-x\otimes 1}.
$$
We denote the residue field of the ideal $\frak m$ by $F$, and let
$D_1$ be the image of the first factor in $F$, and $D_2$ of the
second. It follows from Proposition~\ref{prop_tensor} that the
identity homomorphism $\Qt(D_1)\to \Qt(D)\subseteq L$ under which
$x\otimes 1\mapsto x$ can be extended to an embedding of $F$ into
$L$. We denote the image of an element $1\otimes x$ in $L$ by $x'$.
We can now assume that $D_1$ and $D_2$ are contained in $L$. By
Proposition~\ref{SplitScalExt}, $L$ is the splitting field of the
same family of algebras noth over  $\Qt(D_1)$ and over $\Qt(D_2)$.
By construction of $D_1$ and $D_2$, there exists a differential
isomorphism between them such that $x\mapsto x'$, and therefore also
between their fields of fractions. It follows from
Proposition~\ref{ExtAutoSplit} that this automorphism can be
extended to a differential automorphism of $L$.
\end{proof}

We define the following transformation of families of differentially
finitely generated $K$-algebras $\{B_\alpha\}$. For each ring
$B_\alpha$ consider the family $X_\alpha=\{B_\alpha/\frak m\}$,
where the $\frak m$ are all possible locally maximal differential
ideals of the ring $B_\alpha$. We consider the family
$X=\bigcap_\alpha X_\alpha$; it consists of locally simple
differential algebras. Now for each $B\in X$ there exists $s\in B$
such that $B_s$ is simple. Let $Y$ consist of all the $B_s$ of this
form.

\begin{proposition}
Suppose that $\{\,B_\alpha\,\}$ is an arbitrary family of
differentially finitely generated $K$-algebra, and families $X$ and
$Y$ are obtained in the way described above. Then the splitting
fields of all the three families naturally coincide.
\end{proposition}
\begin{proof}
Let $L$ be the splitting field of the initial family. Consider an
arbitrary differential homomorphism $B_\alpha\to L$. Then the image
of $B_\alpha$ is a locally simple ring (this follows from
Proposition~\ref{constatomic} and the constructedness of the
splitting field) and has the form $B_\alpha/\frak m$ for some
locally maximal differential ideal $\frak m$. On the other hand, if
we are given a differential homomorphism from $B_\alpha/\frak m$
into $L$, then it results in a through homomorphism from $B_\alpha$
to $L$. Therefore if we replace the initial family by the family
$X$, then $L$ will satisfy the definition of a splitting field of
the family $X$. We now point out that after passing to the ring
$B_\alpha/\frak m$ the differential homomorphism thus constructed is
injective. Then we can invert the  corresponding element so that
$(B_\alpha/\frak m)_s$ becomes a simple ring, and the injectivity
implies that the differential homomorphism can be extended from
$B_\alpha/\frak m$ to the entire ring $(B_\alpha/\frak m)_s$.
Therefore in this case the splitting fields also coincide.
\end{proof}

\begin{definition}\label{def_split_pol}
Consider some ring of differential polynomials
$K\{y_1,\ldots,y_n\}$, and let $f_1,\ldots,f_k$ be differential
polynomials in this ring. Then the {\it splitting field of the
system} of differential equations
$$
\left\{
\begin{aligned}
f_1(y_1,&\ldots,y_n)=0\\
&\ldots\\
f_k(y_1,&\ldots,y_n)=0
\end{aligned}
\right.
$$
is defined to be the splitting field of the algebra
$$
B=K\{y_1,\ldots,y_n\}/[f_1,\ldots,f_k].
$$
\end{definition}

We now show that our notation of the splitting field agrees with all
the existing ones.

{\it Splitting field of a polynomial}. If $f$ is a polynomial of one
variable (that is, an element of $K[x]$), then its splitting field
in the usual sense coincides with our splitting field. The
derivations are considered to be zero.

Indeed, let $L$ be the splitting field of the polynomial $f$. We now
show that property~1) in Definition~\ref{def_split} of a splitting
field holds. Consider the ring
$$
R=L\otimes K[x]/(f).
$$
Let $\frak m$ be some locally maximal ideal of this ring. In fact,
it is automatically maximal (see~\cite{AM}, Ch.~5, Exercise~24), but
we can afford not to use this. Let $s$ be an element such that the
ring
$$
F=(R/\frak m)_s
$$
simple, that is, a field. Then $F=L(x)$ by construction, where $x$
is a root of the polynomial $f$. But the fact that $L$ was already
the splitting field implies that $x\in L$ and therefore also $F=L$.
Condition~2) on a splitting field means that it is generated by
roots of $f$; however, this is the case by the definition of the
splitting field of a polynomial. And condition~3) on a splitting
field means that if some field contains all the roots of this
polynomial, then it also contains its splitting field.

{\it Picard-Vessiot extension}. The splitting field of a system of
linear differential equations (or of a linear differential
polynomial) coincides with the corresponding Picard-Vessiot
extension (see~\cite{SW}). Here we assume that all differential
rings have one derivation.

Indeed, let $L$ be the Picard-Vessiot extension for some system of
linear differential equations $y'=Ay$ over a differential field $K$.
We denote the set of these equations by $E$. Consider the ring
$$
R = L \otimes K \{y_1,\ldots,y_n\}/[E].
$$
Now let $\frak m$ be some locally maximal differential ideal of the
ring $R$, and let $s$ be an element such that the ring $(R/\frak
m)_s$ is simple. Then it follows from Proposition~7(d)
in~\textsection~10 of~Ch.~III in~\cite{Kl} that the ring thus
obtained does not contain new constants. Moreover, the images of
$y_1,\ldots,y_n$ form a solution of that equation. Then the
definition of a Picard-Vessiot ring (see~\cite{SW}, Ch.~1,
Definition~1.15 and Lemma~1.7) implies that the vector composed of
the $y_i$ is expressed in terms of bases solutions with constant
coefficients. therefore the elements $y_i$ belong to $L$. Thus, we
have shown that property~1) of Definition~\ref{def_split} holds.
That same definition implies that a Picard-Vessiot ring is generated
by solutions of the system of equations, and therefore property~2)
of Definition~\ref{def_split} holds. Proposition~1.22 and~1.20(2) in
Ch.~1 of~\cite{SW} imply that any Picard-Vessiot field is unique up
to isomorphism and coincides with the field of fractions of the
Picard-Vessiot ring. Thus, if some field contains a fundamental
matrix of solutions of a linear differential equation, then it also
contains the corresponding Picard-Vessiot field.

{\it Parametrized Picard-Vessiot extension}. The splitting field of
a system of linear differential equations (of linear differential
equation) with parameters coincides with the notion of a
parametrized Picard-Vessiot extension (see~\cite{SC}).

Let
$$
\partial_1 Y = A_1 Y,\quad \partial_1 Y = A_2 Y, \quad \ldots,\quad
\partial_r Y = A_r Y
$$
be some systems of equations over the field $K$ with derivations
$\Delta=\{\partial_1,\ldots,\partial_m\}$, where $k\leqslant m$. Let
$L$ be the corresponding parametrized Picard-Vessiot extension. We
denote all equations in the systems above by $E$. Then we consider
the ring
$$
R= L \otimes K\{y_1,\ldots,y_n\}/[E].
$$
Let $\frak m$ be some locally maximal differential ideal of the ring
$R$, and let an element $s\in R$ be such that the algebra $(R/\frak
m)_s$ is simple. Then it follows from Lemma~9.3 in~\cite{SC} that in
the ring $(R/\frak m)_s$ the subring of constants with respect to
the derivations $\partial_1,\ldots,\partial_k$ coincides with the
subfield of constants of the field $K$ with respect to the same
derivations. Therefore by Definition~9.4 in~\cite{SC} this ring is
the corresponding parametrized Picard-Vessiot ring. Therefore, its
field of fractions is the corresponding parametrized Picard-Vessiot
extension. Then a solution composed of the images of the $y_i$ is
expressed in terms of a fundamental matrix with elements in the
field $L$ with coefficients that automatically belong to the field
$K$. Therefore $(R/\frak m)_s$ coincides with $L$. Property~2) in
Definition~\ref{def_split} of a splitting field holds by the
definition of a parametrized Picard-Vessiot extension
(see~\cite{SC}, Definition~9.4(2b)). Property~3) in
Definition~\ref{def_split} of a splitting field means that any
differential field generated by a fundamental matrix of solutions of
the equations in $E$ is isomorphic to the parametrized
Picard-Vessiot extension; but this is Proposition~9.5(1)
in~\cite{SC}.

{\it Strongly normal extensions}. Strongly normal extensions in the
sense of Kolchin are also splitting fields of the corresponding
differential algebras (see~\cite{JK}).

Let $G$ be a strongly normal extension of a differential field $K$
in the sense of Definition~12.1 in~\cite{JK}. Then by
Proposition~13.8 in~\cite{JK} there exists a simple differential
algebra $R\subseteq G$ finitely generated over $K$ such that $G$ is
the field of fractions of $R$. We claim that $G$ is a splitting
field of the algebra $R$. Let us show that property~1) of
Definition~\ref{def_split} of a splitting field holds. Consider the
algebra $G\otimes R$. We can pass to its localization $P=G\otimes G$
without loss of prime differential ideals (see Corollary~13.6
in~\cite{JK}). It now follows from Corollary~14.3 in~\cite{JK} that
any locally maximal differential ideal $\frak m$ of $P$ is maximal;
moreover, in the proof of that corollary it was shown that $P/\frak
m=G$. Property~2) of Definition~\ref{def_split} of a splitting field
holds, since $G$ is the field of fractions of $R$. If some field
contains the algebra $R$, then it automatically contains also its
field of fractions, that is, the field $G$, and this precisely means
that property~3) of Definition~\ref{def_split} of a splitting field
holds.

\subsection{Splitting subfields}\label{sec72}

Let $\{B_\alpha\}$ be some family of differentially finitely
generated $K$-algebras. Recall that we can consider all possible
differential homomorphisms from the rings $B_\alpha$ into some
differential field $F$. For convenience we denote the images of such
homomorphisms by $\overline{B}_\alpha$. Let $\overline{K}$ be the
differential closure of the field $K$, and $L$ some splitting field
over $K$. Then property~3) of Definition~\ref{def_split} of a
splitting field guarantees the existence of an embedding of $L$ into
$\overline{K}$. However, generally speaking, there may exist many
such embeddings; in this connection we highlight the following
proposition.

\begin{proposition}\label{SplitSub}
Let $L$ be the splitting field over $K$ of a family of
differentially finitely generated $K$-algebras $\{B_\alpha\}$. Then
among the differential subfields of $\overline{K}$ isomorphic to $L$
there exists the largest one, and it coincides with the smallest
subfield generated by rings of the form $\overline{B}_\alpha$.
\end{proposition}
\begin{proof}

Consider the subfield of $\overline{K}$ generated by all
differential subrings isomorphic to $\overline{B}_\alpha$. We denote
this field by $L'$. If this subfield is the splitting field, then it
is the required one. We know from Proposition~\ref{ExistSplitClose}
that the field $\overline{L}$ can be mapped isomorphically onto
$\overline{K}$. We claim that $L$ is mapped to $L'$ under this
isomorphism. Since by definition $L$ is generated by the
differential algebras isomorphism to $\overline{B}_\alpha$, its
image is contained in $L'$. Conversely, for any algebra $B$ of the
form $\overline{B}_\alpha$ it is true that $L\cdot B$ is locally
simple, and therefore by definition of $B$ is contained in $L$.
\end{proof}

\begin{definition}\label{def47}
The largest subfield of $\overline{K}$ isomorphic to some splitting
field over $K$ is called a {\it splitting field}.
\end{definition}

\begin{definition}\label{def48}
We call a subfield $L$ of the field $\overline{K}$ containing $K$
{\it good} if $\overline{K}$ is the differential closure of the
field $L$. By Corollary~\ref{DCFuniversUnique} a subfield $L$ of the
field $\overline{K}$ is {\it good} if and only if the field
$\overline{K}$ is constructed over $L$.
\end{definition}

We point out that in the case of a Picard-Vessiot extension and a
strongly normal extension of fields, any differential subfield is
good. Indeed, each of these extensions is finitely generated, and,
consequently, any differential subfield is also differentially
finitely generated. The Proposition~\ref{FiniteConstr} guarantees
that all of them are good. In the case of Parametrized
Picard-Vessiot extensions Theorem~3.5(3) in~\cite{SC} ensures that
all differential subfields are good. It should be noted that in the
general case not all subfields are good. However, the following
proposition holds.

\begin{proposition}\label{SubSplitConstr}
If $L$ is a splitting subfield, then it is good.
\end{proposition}
\begin{proof}
Since the differential closure of $L$ coincides with the
differential closure of $K$ (Proposition~\ref{ExistSplitClose}), we
obtain what is required.
\end{proof}

We define the {\it full differential Galois group} $\diffgal
(\overline{K}/K)$ of the field $K$ to be the group of all
differential automorphisms of $\overline{K}$ over $K$. Since the
filed $K$ is fixed everywhere, we denote this group by $G$. We say
that a differential subfield $L$ is invariant in $\overline{K}$ over
$K$ if $G(L)\subseteq L$.

\begin{proposition}\label{NormSplitCrit}
A differential subfield $L$ is invariant in $\overline{K}$ over $K$
if and only if it is the splitting subfield of some family.
\end{proposition}
\begin{proof}
A subfield is invariant if and only if, together with any simple
differentially finitely generated algebra, this subfield contains
all algebras isomorphism to this algebra. But
Proposition~\ref{SplitSub} implies that these are precisely
splitting subfields.
\end{proof}

\begin{corollary}
If $L$ is differential subfield invariant in $\overline{K}$ over
$K$, then $\overline{K}$ is constructed over $L$.
\end{corollary}
\begin{proof}

If the field $L$ is invariant in $\overline{K}$, then by
Proposition~\ref{NormSplitCrit} it is a splitting subfield of
$\overline{K}$. Then Proposition~\ref{SubSplitConstr} completes the
proof.
\end{proof}

\subsection{Normal extensions}\label{sec73}

In this subsection we deal with differential fields contained not in
the differential closure but in an arbitrary splitting field. We say
that an extension of differential fields $F\subseteq L$ is {\it
normal} if $L$ is the splitting field of some family of
differentially finitely generated algebras over $F$.

\begin{proposition}\label{prop_52}
An extension of differential fields $F\subseteq L$ is normal if and
only if $L$ can be embedded into $\overline{F}$ as an differential
subfield that is invariant over $F$.
\end{proposition}
\begin{proof}

Suppose that $F\subseteq L$ is normal; then by
Proposition~\ref{ExistSplitClose} the differential closure of $L$
coincides with $\overline{F}$. We claim that $L$ is invariant in it.
Let $L$ be the splitting field of a family $\{\,B_\alpha\,\}$, and
let $B$ be a differential subalgebra of $\overline{F}$ isomorphic to
some $\overline{B}_\alpha$. We know that the differential closure is
constructed over the field $L$. Then by
Proposition~\ref{constatomic} the algebra $L\cdot B$ is locally
simple and, consequently, coincides with $L$. Consequently, $L$ is a
splitting subfield of $\overline{F}$, which is invariant by
Proposition~\ref{NormSplitCrit}.

Conversely, let $L$ be embedded as an invariant subfield into
$\overline{F}$. Then it follows from Proposition~\ref{NormSplitCrit}
that $L$ is a splitting subfield and, a fortiori, a splitting field;
consequently, $L$ is normal.
\end{proof}

It is important that for normal extensions all the same propositions
hold as for the differential closure. Then it is not necessary to
work with the entire differential closure as a whole; we can confine
ourselves to the normal extension we are interested in.

\begin{proposition}\label{NormalSub}
Let $L$ be some splitting field over $K$, and let $B$ be a simple
differentially finitely generated $K$-algebra contained in $L$. Then
there exists a map from the splitting field $F$ of the ring $B$ into
$L$. Among the differential subfields of $L$ isomorphic to $F$ there
exists the largest one, and it coincides with the smallest subfield
generated by all the differential rings isomorphic to $B$.
\end{proposition}
\begin{proof}

We need to map the field $F$ into $L$. We embed $L$ into its
differential closure as a splitting subfield over the ground field
$K$. Then $L $ is invariant under the action of the group $G$, that
is, together with the algebra $B$, it also contains any algebra
isomorphic to $B$ (Proposition~\ref{SplitSub}). In other words, $L$
contains the smallest subfield spanned by the differential algebras
isomorphic to $B$. By Proposition~\ref{SplitSub} the latter subfield
is isomorphic to $F$. By construction it is the largest one.
\end{proof}

\begin{definition}
If $L$ is a splitting field over $K$, and $F$ the largest
differential subfield of $L$ isomorphic to some splitting field over
$K$, then $F$ is called a {\it splitting subfield}.
\end{definition}

From Propositions~\ref{prop_52},~\ref{NormalSub} we obtain a simple
corollary.

\begin{corollary}\label{SplitNormSplit}
Every splitting field is normal over its a splitting subfield.
\end{corollary}

\begin{proposition}\label{stablenorm}
Let $L$ be arbitrary splitting field over $K$, and $F$ an arbitrary
splitting subfield over $K$. Then $F$ is invariant under any
differential automorphism of $L$ over $K$.
\end{proposition}
\begin{proof}
The proof immediately follows from Proposition~\ref{NormalSub}.
\end{proof}

The {\it differential Galois group} $\diffgal(L/F)$ of a normal
extension $F\subseteq L$ is defined to be the group of all
differential automorphisms of $L$ over $F$.

\begin{proposition}\label{InvarSubSplitCrit}
A differential subfield $F$ of a splitting field $L$ over $K$ is
invariant under $\diffgal(L/K)$  if and only if it is the splitting
subfield of some family of differentially finitely generated
$K$-algebras.
\end{proposition}
\begin{proof}
We embed $L$ into its differential closure $\overline{L}$. It
follows from Proposition~\ref{ExtAutoSplit} that any differential
automorphism of $L$ can be extended to a differential automorphism
of $\overline{L}$. Therefore the subfield $F$ is invariant in $L$ if
and only if it is invariant in $\overline{L}$. But by
Proposition~\ref{NormSplitCrit} this is equivalent to $F$ being a
splitting subfield of $\overline{L}$. Then it follows from
Proposition~\ref{NormalSub} that if some simple differentially
finitely generated subalgebra $B$ is contained in the field $L$,
then $L$ contains all differential subalgebras of $\overline{L}$
isomorphic to $B$. Therefore a subfield of $L$ is a splitting
subfield if and only if it is a splitting subfield of
$\overline{L}$.
\end{proof}

\begin{proposition}\label{NormInvar}
Let $F\subseteq L$ be some normal extension of differential fields.
Then $F=L^{\diffgal (L/F)}$.
\end{proposition}
\begin{proof}
For proving this we need to apply Proposition~\ref{InvariantSplit}
with $F$ instead of $K$.
\end{proof}

\section{Galois correspondence for normal extensions}\label{sec8}

We must point out that we construct the theory under the assumption
that all differential rings are algebras over some differential
field $L$ of characteristic zero. As a natural such field one can
always choose the field of rational numbers $\mathbb Q$.

\begin{theorem}\label{GaloisCorr}
Let $K\subseteq L$ be some normal extension of differential fields.
Then for any splitting subfield $F$ of $L$ containing  $K$ we have
the equation
$$
F = L^{\diffgal(L/F)}.
$$
Moreover, for the group of differential automorphisms we have the
equation
$$
\diffgal(F/K)=\diffgal(L/K)/\diffgal(L/F).
$$
\end{theorem}
\begin{proof}
Corollary~\ref{SplitNormSplit} says that the extension $F\subseteq
L$ is normal. Then by Proposition~\ref{NormInvar} we obtain that the
field $F$ can be reconstructed from its differential Galois group.

Proposition~\ref{stablenorm} guarantees that any differential
automorphism in the group $\diffgal(L/K)$ is correctly restricted to
the subfield $F$. Then the kernel of this homomorphism is exactly
the group $\diffgal(L/F)$. It remains to observe that any
differential automorphism of the field $F$ can be extended to a
differential automorphism of the field $L$. Indeed, this follows
from Proposition~\ref{ExtAutoSplit}.
\end{proof}

We now introduce the following notation. Let $F\subseteq L$ be some
normal extension of differential fields; in other words, $L$ is some
splitting field over $F$. We denote by $\mathcal F$ the set of all
good subfields in $L$ containing $F$. We denote by $\mathcal F'$ the
set $\{\,L^H\mid H\subseteq \diffgal (L/F)\,\}$. Let $\mathcal N$
denote the set of subfields invariant $\diffgal(L/F)$.

Note that the set $\mathcal F'$, generally speaking, does not
contain all the intermediate subfields between $F$ and $L$. The
definition of a good subfield implies the inclusion  $\mathcal
F\subseteq \mathcal F'$, which, generally speaking, is strict. The
inclusion $\mathcal N\subseteq\mathcal F$ is also always true (see
Proposition~\ref{SubSplitConstr}).

\begin{example}
This example is borrowed from~\cite{Rl} and~\cite{KL3}. We show how
bad a normal extension can be, and together with it its Galois
group.

We confine ourselves to the case of one derivation. Let the field
$K$ be an algebraically closed field of constants. We consider the
equation $y'=y^3-y^2$. Then it follows from Corollary on p.~532
in~\cite{Rl} that its splitting field is the following field:
$$
L=K(x_1,\ldots,x_n,\ldots),
$$
where $x_k$ are algebraically independent over $K$ and
$x_k'=x^3_k-x^2_k$. As we see, its differential Galois group is the
group of all permutations of the elements, that is, the permutation
group $S_{\mathbb N}$ of the set of positive integers.

The subfields generated by finitely many $x_i$, that is fields of
the form $K(x_{i_1},\ldots,x_{i_n})$, are good and are contained in
$\mathcal F$. The subfield $K(x_{2},x_4,\ldots,x_{2n},\ldots)$ is
contained in $\mathcal F'\setminus\mathcal F$. Indeed, it is the set
of fixed points for the subgroup permuting only the $x_k$ with odd
indices. But this subfield itself is already a splitting field, and,
consequently, $L$ cannot  be constructed over it. The subfield
$K(x_2,x_3,\ldots,x_n,\ldots)$ is not contained even in $\mathcal
F'$, since any permutation that is the identity on all elements
except the first one is the identity everywhere.

We define a subgroup  $A^F_{\mathbb N}$ of $S_{\mathbb N}$ as
follows: a permutation $\sigma$ belongs to $A^F_{\mathbb N}$ if an
only if it permutes only finitely many elements by an even
permutation. This group is the smallest normal subgroup of
$S_{\mathbb N}$. However its field of invariants is equal to $K$.
Therefore, in $L$ there is not any proper splitting subfield. In
other words, for any element $f\in L$, by using the operations of
addition, subtraction, multiplication, division, derivations, and
changing variables $x_i$l to $x_j$, one can obtain any other element
of $L$.
\end{example}

\section{Connection with differential algebraic
varieties}\label{sec9}

Let $B$ be an arbitrary differentially finitely generated algebra
over the field $K$. Let $\overline{K}$ be the differential closure
of the field $K$. We denote by $X_{\overline{K}}$ the set of
$\overline{K}$points of the algebra $B$, that is, the set of
differential homomorphisms from $B$ into $\overline{K}$ over $K$. As
above, we denote by $G$ the group $\diffgal(\overline{K}/K)$. The
group $G$ naturally acts on $X_{\overline{K}}$:
$$
\xi\mapsto g\circ\xi,\quad g\in G,\quad \xi\in X_{\overline{K}}.
$$
We indicate the following connection of $\diffsmax B$ and
$X_{\overline{K}}$.

\begin{theorem}\label{Variety}
Under the conditions stated above, we have the equality
$$
\diffsmax B=X_{\overline{K}}/G.
$$
\end{theorem}
\begin{proof}

We construct a map from $X_{\overline{K}}$ into $\diffsmax B$;
namely, with each $\xi$ we associate its kernel $\ker\xi$. This map
is well-defined, since by Corollary~\ref{DCFuniversUnique} the field
$\overline{K}$ is constructed over $K$ and therefore by
Proposition~\ref{constatomic} the image of $B$ is locally simple.
Since any locally simple algebra embeds into the differential
closure, the map thus constructed is surjective. By the definition
of the action of the group $G$ this map is constant on the orbits of
the action of the group; consequently, we have a map
$$
X_{\overline{K}}/G\to\diffsmax B.
$$
It remains to show that it is injective. Indeed, suppose that two
differential homomorphisms $\xi_1,\xi_2\colon B\to \overline{K}$
have the same kernels. We denote by $B_1$ and $B_2$ the images of
$B$ under the action of $\xi_1$ and $\xi_2$, respectively:
$$
\xymatrix@R=6pt{
    B_1\ar@{}[r]|{\subseteq}&\overline{K}&B_2\ar@{}[l]|{\supseteq}&\\
    &B\ar[ul]^{\xi_1}\ar[ur]_{\xi_2}&\\
}
$$
then the differential algebra $B_1$ is isomorphic to $B_2$, and this
isomorphism extends to their fields of fractions $F_1$ and $F_2$,
respectively. By Proposition~\ref{FiniteConstr} the field
$\overline{K}$ is constructed over $F_1$ and $F_2$ and coincides
with their differential closure. Consequently, by
Corollary~\ref{DCFuniversUnique} of the uniqueness theorem this
isomorphism can be extended to a differential automorphism $g$ of
the field $\overline{K}$. By construction it satisfies the property
$\xi_2=g\circ\xi_1$, which completes the proof.
\end{proof}

Theorem~\ref{Variety} indicates a connection between a differential
algebraic variety and the set of locally closed points of the
differential spectrum. In the study of algebraic varieties by the
methods of scheme theory it is very important that the properties of
the spectrum are in a sense related to the properties of the
variety. This theorem asserts that a similar relation also exists in
the case of differential algebraic varieties. Thus, it is possible
to transfer geometric results on the differential spectrum to the
case of differential algebraic varieties.

\end{document}